\documentclass[a4paper,12pt,reqno,twoside]{amsart}

\usepackage[utf8]{inputenc} 		
\usepackage[T1]{fontenc} 

\usepackage{amssymb}
\usepackage{amsmath}
\usepackage{graphicx}
\usepackage[colorlinks,citecolor=blue,urlcolor=blue]{hyperref}
\usepackage{cite}
\usepackage[hmargin=2.5cm,vmargin=2.5cm]{geometry}
\usepackage{mathrsfs} 
\usepackage{ellipsis}

\newcommand{\be}{\begin{eqnarray}}
\newcommand{\ee}{\end{eqnarray}}

\newcommand{\beq}{\begin{equation}}
\newcommand{\eeq}{\end{equation}}
\newcommand{\beqn}{\begin{equation*}}
\newcommand{\eeqn}{\end{equation*}}

\newcommand{\slot}{\,\cdot\,}

\DeclareMathOperator{\tr}{tr}

\newtheorem{thm}{Theorem}[section]
\newtheorem{prop}[thm]{Proposition}
\newtheorem{cor}[thm]{Corollary}
\newtheorem{lem}[thm]{Lemma}
\newtheorem{defn}[thm]{Definition}
\newtheorem{sublem}[thm]{Sublemma}

\newcommand\cB{{\mathcal B}}

\newcommand\cN{{\mathcal N}}

\newcommand\cQ{{\mathcal Q}}

\newcommand\cS{{\mathcal S}}
\newcommand\cT{{\mathcal T}}

\newcommand\bE{{\mathbb E}}

\newcommand\bN{{\mathbb N}}

\newcommand\bR{{\mathbb R}}

\newcommand\rd{{\mathrm d}}

\newcommand{\ve}{\varepsilon}

\begin{document}

\title{Stein's method for dynamical systems}

\author[Olli Hella]{Olli Hella}
\address[Olli Hella]{
Department of Mathematics and Statistics, P.O.\ Box 68, Fin-00014 University of Helsinki, Finland.}
\email{olli.hella@helsinki.fi}

\author[Juho Lepp\"anen]{Juho Lepp\"anen}
\address[Juho Lepp\"anen]{
Department of Mathematics and Statistics, P.O.\ Box 68, Fin-00014 University of Helsinki, Finland.}
\email{juho.leppanen@helsinki.fi}

\author[Mikko Stenlund]{Mikko Stenlund}
\address[Mikko Stenlund]{
Department of Mathematics and Statistics, P.O.\ Box 68, Fin-00014 University of Helsinki, Finland.}
\email{mikko.stenlund@helsinki.fi}
\urladdr{http://www.helsinki.fi/~stenlund/}


 





\begin{abstract}
We present an adaptation of Stein's method of normal approximation to the study of both discrete- and continuous-time dynamical systems. We obtain new correlation-decay conditions on dynamical systems for a multivariate central limit theorem augmented by a rate of convergence. We then present a scheme for checking these conditions in actual examples. The principal contribution of our paper is the method, which yields a convergence rate essentially with the same amount of work as the central limit theorem, together with a multiplicative constant that can be computed directly from the assumptions.
\end{abstract}

\maketitle


\subsection*{Acknowledgements}
This work was supported by the Jane and Aatos Erkko Foundation, and by Emil Aaltosen S\"a\"ati\"o. 


%
%
%
%
%
%

\section{Introduction}\label{sec:intro}

Let $T$ be a measure-preserving transformation on a probability space $(X,\cB,\mu)$, and let $f:X\to\bR^d$ be a function with $d\ge 1$ and $\int f\,\rd\mu = 0$. Then the iterates of $T$ determine a sequence of centered random vectors $f\circ T^k$, $k\ge 0$, the randomness arising solely from the distribution~$\mu$ on the sample space~$X$. It is common to interpret~$f$ as an observable quantity,~$T$ as the rule of deterministic time evolution, and $f\circ T^k$ as the value of the observable at time~$k$. This paper concerns the problem of approximating the law of the scaled time averages
\beqn
W(N) = \frac{1}{\sqrt N}\sum_{k=0}^{N-1} f\circ T^k
\eeqn
by a normal distribution of $d$ variables: establish conditions on~$T$ which imply an upper bound on the distance between the two distributions for a given value of~$N\ge 1$, with respect to a suitable metric. (For simplicity, we restrict to discrete time in this section.)

In the present context of dynamical systems (i.e., measure-preserving transformations), the question has, of course, been long studied in a number of landmark papers, e.g., \cite{Jan_2000,Pene_2002,Pene_2004,Pene_2005,Gouezel_2005,Dubois_2011,GouezelMelbourne_2014}, using different methods, hypotheses, and metrics. To our knowledge, however, P\`ene's important article~\cite{Pene_2005} is the only one treating the multidimensional case~$d>1$, building on a method due to Rio~\cite{Rio_1996}. We note that while a multivariate central limit theorem can be reduced to the univariate case, bounds on the rate of convergence do not pass from~$d=1$ to higher dimensions.

In the seminal paper~\cite{Stein_1972}, Stein introduced a new and unusual method for (univariate) normal approximation, which has in recent years received a great amount of attention in the probability theory literature. This paper grew out of interest in investigating how it can be adapted to the setup of dynamical systems. Stein's method stands out among other methods of normal approximation known to us at least in two regards:\@ Firstly, it does not resort to the use of characteristic functions (Fourier transforms) of the distributions, even under the surface (in contrast with, say, the martingale central limit theorem or the method used in~\cite{Pene_2005}). Rather, it allows us to turn the problem of normal convergence into a set of problems on correlation decay. Secondly, on top of establishing normal convergence, Stein's method produces a rate of convergence for free. We submit that both of these features are useful in the study of dynamical systems. In essence, additionally to connecting Stein's method to dynamical systems, which we believe is valuable in its own right, our paper is about bridging the gap between proving central limit theorems and obtaining rates of convergence, all based on information on correlation decay.

We mention a further benefit of the above features of Stein's method: the method turns out to be sufficiently flexible to be implemented in the non-stationary setup in which the random vectors $f\circ T^k$ are replaced by $f\circ T_k\circ\cdots \circ T_1$ with a rather arbitrary sequence of maps~$T_i$, $i\ge 1$, which generally do not preserve a common measure. This appears to be a big difference compared to P\`ene's adaptation \cite{Pene_2005} of Rio's method. (Further differences will be discussed in Section~\ref{sec:billiards}.) Such a substantial generalization is currently work in progress, and will be reported on later as part of Olli Hella's doctoral thesis. It will in particular apply to so-called quasistatic dynamical systems, recently studied in~\cite{DobbsStenlund_2016,Stenlund_2016,LeppanenStenlund_2016}.

Let us pause and briefly point out a curious aspect of merging Stein's method and dynamical systems. It is well known to experts on the statistical properties of dynamical systems that adding random noise to the rules of the dynamics facilitates proving limit theorems.
On the other hand, in its original setting of purely random processes, Stein's method can be challenging to apply, as it involves elements of both ingenuity and art. For example, a typical step required is that of ``auxiliary randomization'', which calls for an exchangeable-pair or a (zero-bias, size-bias, etc.\@) coupling construction, but this works only under favorable circumstances. One is then perhaps tempted to doubt the applicability of the method in the context of deterministic dynamics. Nevertheless, for want of a more precise expression, there is a certain flavor of rigidity to deterministic systems, which can be taken advantage of. It allows, at least for systems of hyperbolic type, to circumvent auxiliary randomization and resort to a direct approach. We refer to the main theorems and Section~\ref{sec:applications} for a more comprehensive picture.

Our paper appears to be the first general attempt to adapt Stein's method of normal approximation to the context of dynamical systems. In~\cite{Gordin_1990}, an abstract univariate central limit theorem is obtained, via an idea due to Stein, for automorphisms $T$ admitting a so-called homoclinic transformation. The latter theory is known to apply in a few examples: hyperbolic toral automorphisms and Bernoulli shifts~\cite{Gordin_1990}, and a zero entropy transformation~\cite{King_1997}. Moreover, the application of Stein's method in~\cite{Gordin_1990,King_1997} is only partial, in that --- quoting~\cite{Gordin_1990} --- it does \emph{``not derive any estimates of the rate of convergence in the CLT''}. Stein's method has been implemented in the different dynamical systems setup of Poisson limits, in connection with hitting time statistics for shrinking sets; see~\cite{Denker_etal_2004, Psiloyenis_2008, GordinDenker_2012, Haydn_2013, HaydnYang_2016}.

Finally, we remark that the latest probability theory literature on Stein's method has to a large extent focused on multivariate normal approximation; see, e.g.,~\cite{Gotze_1991,RinottRotar_1996,GoldsteinRinott_1996,ChatterjeeMeckes_2008,Meckes_2009,ReinertRollin_2009,Nourdin_etal_2010,Chen_etal_2011,Meckes_2012,Rollin_2013,Gaunt_2016,Berckmoes_etal_2016}. Our adaptation of Stein's method covers the multivariate case.

\smallskip
\noindent{\bf Notations and conventions.} Throughout this paper, $T:X\to X$ is a measure-preserving transformation on a probability space \((X,\cB,\mu)\). Given a function \(f: \, X \to \bR^d\), we denote its mean by \(\mu(f) = \int_{X} f \, d\mu\). The coordinate functions of \(f\) are denoted by \(f_{\alpha}\), $\alpha \in\{1,\dots,d\}$. We define
\beqn
\|f\|_\infty = \max_{1\le\alpha\le d}\|f_\alpha\|_\infty.
\eeqn
We also set \(f^n = f \circ T^n\) for all \(n \in \bN_0 = \{0,1,2,\ldots\}\), with \(f^0 = f\).

For a function \(B: \, \bR^d \to \bR^{d'}\), we write \(D^kB\) for the $k$th derivative. We define
\begin{align*}
\Vert D^k B \Vert_{\infty} = \max \{ \|\partial_1^{t_1}\cdots\partial_d^{t_d}B_\alpha\|_{\infty}  :  t_1+\cdots+t_d = k,\, 1\le \alpha\le d'  \}.
\end{align*}
Here $B_\alpha$, $1\le \alpha\le d'$ are the coordinate functions of~$B$.
The gradient notation $\nabla B = DB$ is reserved for the scalar-valued case $d'=1$. 

Finally, given two vectors $v,w\in\bR^d$, we write $v\otimes w$ for the $d\times d$ matrix with entries
\beqn
(v\otimes w)_{\alpha\beta} = v_\alpha w_\beta.
\eeqn

\bigskip
\noindent{\bf Structure of the paper.} In Section~\ref{sec:results} we present our main results in both discrete and continuous time. The results are obtained using Stein's method, which is reviewed in Section~\ref{sec:Stein}. Since the cases of discrete and continuous time are virtually identical, Sections~\ref{sec:preliminary} and~\ref{sec:proofs} are devoted solely to discrete time: In Section~\ref{sec:preliminary} we formulate preliminary versions of the main results; the former are proved in the same section using the method of Section~\ref{sec:Stein}. The main results are then inferred from these in Section~\ref{sec:proofs}. Then, in Section~\ref{sec:proofs_flow}, we prove the continuous-time results. Finally, in Section~\ref{sec:applications} we present an abstract scheme for verifying a key assumption of our theorems. We proceed to give two examples, the first one of which is trivial --- but already reveals the essential ideas --- while the second one is more technical in nature. Although deferred to the very end, the scheme of Section~\ref{sec:applications} constitutes an integral part of our work, as it indicates the usage of the theorems.


\section{Main results}\label{sec:results}

\subsection{Discrete time: multivariate case}
In order to state our main results, we introduce some definitions. 

Let \(f : \, X \to \bR^d\) be given. Given \(N\in \bN_0\), we write 
\beqn
W = W(N) = \frac{1}{\sqrt{N}}\sum_{k=0}^{N-1} f^k.
\eeqn
Given also $K\in \bN_0\cap [0,N-1]$, we write
\beqn
[n]_K = [n]_K(N) = \{ k \in \bN_0 \cap [0,N-1] \, : \, |k-n| \le K \}
\eeqn
and
\beq\label{eq:W^n}
W^n = W^n(N,K) = W - \frac{1}{\sqrt{N}} \sum_{k\in [n]_K} f^k
\eeq
for all $n\in\bN_0\cap [0,N-1]$. In other words, $W^n$ differs from~$W$ by a time gap (within $[0,N-1]$) of radius~$K$, centered at time~$n$. Observe that the cardinality of $[n]_K$ satisfies
\beqn
K+1 \le |[n]_K| \le 2K + 1.
\eeqn

Given a unit vector $v\in\bR^d$, we say that $f$ is a \emph{coboundary in the direction~$v$} if there exists a function $g_v:X\to\bR$ in $L^2(\mu)$ such that
\beqn
v\cdot f = g_v - g_v\circ T.
\eeqn
In the scalar case $d=1$ we simply say that $f$ is a coboundary.

We denote by $\Phi_\Sigma(h)$ the expectation of a function $h:\bR^d\to\bR$ with respect to the $d$-dimensional centered normal distribution $\cN(0,\Sigma)$ with covariance matrix~$\Sigma\in\bR^{d\times d}$:
\beqn
\Phi_\Sigma(h) = \frac{1}{\sqrt{2\pi\det\Sigma}} \int_{\bR^d} e^{-\frac12 w\cdot\Sigma^{-1}w}h(w)\, dw.
\eeqn

Our first theorem concerns $\bR^d$-valued functions~$f$ and multivariate normal approximation of~$W$ with arbitrary $d\ge 1$:

\begin{thm}\label{thm:main} Let \(f: \, X \to \bR^d\) be a bounded measurable function with \(\mu(f) = 0\). Let \(h: \, \bR^d \to \bR\) be three times differentiable with \(\Vert D^k h \Vert_{\infty} < \infty\) for \(1 \le k \le 3\). Fix integers $N>0$ and $0\le K<N$. Suppose that the following conditions are satisfied:

\begin{itemize}
\item[\bf (A1)] There exist constants $C_2 > 0$ and $C_4 > 0$, and a non-increasing function \(\rho : \, \bN_0 \to \bR_+\) with \(\rho(0) = 1\) and \( \sum_{i=1}^{\infty} i\rho(i) < \infty\), such that
\begin{align*}
|\mu(f_{\alpha} f^k_{\beta})| & \le C_2\,\rho(k)
\\
|\mu(f_{\alpha} f_{\beta}^l f_{\gamma}^m f_{\delta}^n)| & \le C_4 \min\{\rho(l),\rho(n-m)\}
\\
|\mu(f_{\alpha} f_{\beta}^l f_{\gamma}^m f_{\delta}^n) - \mu(f_{\alpha} f_{\beta}^l)\mu(f_{\gamma}^m f_{\delta}^n)| & \le C_4\,\rho(m-l)
\end{align*}
hold whenever $k\ge 0$; $0\le l\le m \le n < N$; $\alpha,\beta,\gamma,\delta\in\{\alpha',\beta'\}$ and $\alpha',\beta'\in\{1,\dots,d\}$.

\smallskip
\item[\bf (A2)] There exists a function \(\tilde\rho : \, \bN_0 \to \bR_+\) such that
\begin{align*}
|\mu( f^n \cdot  \nabla h(v + W^n t) ) | \le \tilde\rho(K)
\end{align*}
holds for all $0\le n < N$, $0\le t\le 1$ and $v\in\bR^d$.

\smallskip
\item[\bf (A3)] $f$ is not a coboundary in any direction. 

\end{itemize}
Then
\beq\label{eq:Sigma}
\Sigma = \mu(f \otimes f) + \sum_{n=1}^{\infty} (\mu(f^n \otimes f) + \mu(f \otimes f^n))
\eeq
is a well-defined, symmetric, positive-definite, $d\times d$ matrix; and
\beq\label{eq:main}
|\mu(h(W)) - \Phi_{\Sigma}(h)| \le C_* \! \left(\frac {K+1}{\sqrt{N}} + \sum_{i= K+1}^{\infty}  \rho(i)\right)  + \sqrt N \tilde\rho(K),
\eeq
where
\beq\label{eq:main_constant}
C_* = 12d^3\max\{C_2,\sqrt{C_4}\}\left(\Vert D^2 h\Vert_{\infty} + \Vert f \Vert_{\infty} \Vert D^3 h \Vert_{\infty} \right)\sum_{i = 0}^\infty (i+1)\rho(i) 
\eeq
is independent of $N$ and $K$.
\end{thm}

The theorem is proved in Sections~\ref{sec:preliminary} and~\ref{sec:proofs}. Observe that both the assumptions as well as the conclusion concern only the given functions~$f$ and~$h$, and the fixed numbers~$N$ and~$K$. Moreover, the bound in~\eqref{eq:main} is expressed purely in terms of the quantities appearing in the assumptions. A smaller value for the constant $C_*$ could be achieved, but we have opted for a clean expression instead.

Note that Assumption~(A1) requires decay of correlations of orders two and four, at a rate which has a finite first moment. Observe also that the second order mixing condition in~(A1) and the condition $\sum_{i=1}^\infty i\rho(i)<\infty$ are the only conditions in the theorem required to hold for arbitrarily large times. These facilitate specifying the target distribution $\cN(0,\Sigma)$ and simplifying the error estimate.

The formulation of Assumption~(A2) in terms of the \emph{function}~$\tilde\rho$ is intentional, even though~$K$ is initially assumed fixed. For the bound in~\eqref{eq:main} to be of any use, we need $K\ll \sqrt N$, and $\tilde\rho(K)$ has to be small. The latter requires that the random vectors
\beq\label{eq:nabh}
f^n
\quad\text{and}\quad
 \nabla h(v + W^n t) =  \nabla h\!\left(v + \frac{t}{\sqrt{N}}\sum_{0\le k< n-K} f^k + \frac{t}{\sqrt{N}}\sum_{n+K<k< N} f^k\right)
\eeq
are nearly orthogonal. Since $\mu(f^n) = 0$, this in turn follows if they are componentwise nearly uncorrelated. Inspecting~\eqref{eq:nabh}, it is thus conceivable that $\tilde\rho(K)$ decreases as $K$ increases, provided~$T$ has sufficiently good mixing properties, because $\nabla h(v + W^n t)$ is a smooth bounded function of the random vectors~$f^k$ with~$|k-n|>K$~only. Section~\ref{sec:applications}, which is an important part of our paper, contains a scheme for verying Assumption~(A2), as well as some examples.
In fact, in these examples, $\rho(i)$ and $\tilde\rho(i)$ have similar asymptotic behavior in the limit $i\to\infty$, in line with intuition.

Note also that the left side of \eqref{eq:main} is independent of~$K$, while on the right side there are competing (increasing and decreasing) contributions involving~$K$. These facts reflect how the theorem is applied in practice: one verifies~(A2) for a range of values of~$K$, and then optimizes the choice of~$K$ as a function of~$N$, in order to obtain a good upper bound in~\eqref{eq:main} which only depends on~$N$. We will return to this point shortly, in Corollary~\ref{cor:exponential}.

Assumption (A3) has to do with $\Sigma$ being positive-definite: Suppose $\Sigma$ is only semi-positive-definite. Then there exists a unit vector $v\in\bR^d$ such that $v\cdot\Sigma v = 0$. Defining $f_v = v\cdot f$, note that $0 = v\cdot \Sigma v = \mu(f_v\,f_v) + 2\sum_{n=1}^{\infty} \mu(f_v \,f_v\circ T^n)$. It follows from \cite{Robinson_1960,Leonov_1961} (together with (A1)) that there exists a function $g_v\in L^2(\mu)$ such that $f_v = g_v - g_v\circ T$, i.e., $f$ is a coboundary in the direction~$v$. The latter is a very restrictive condition on~$f$. For instance, in many applications it is known that such a function~$g_v$ must have a regular representative, so $\sum_{k = 0}^{p-1} f_v(T^k x) = 0$ for all periodic points~$x$ of all periods~$p$, and that periodic points are dense in~$X$.

Coming back to the choice of~$K$ depending on~$N$, we make a simple observation:
\begin{cor}\label{cor:exponential}
Let $\rho(i) = \lambda^i$ and $\tilde\rho(i) = \widetilde C\lambda^i$, where $\lambda\in(0,1)$ and $\widetilde C>0$, and suppose the assumptions of Theorem~\ref{thm:main} are satisfied for $K = \lceil{\log N/|\log\lambda|}\rceil$ and some $N > 2$.
Then
\beqn
|\mu(h(W)) - \Phi_{\Sigma}(h)| 
\le\emph{const}\cdot \frac{\log N}{\sqrt N}
\eeqn 
with $\emph{const} = C_* \! \left(\frac{2}{|{\log\lambda}|} + \frac{\lambda}{\sqrt 3(1-\lambda)}\right) + \widetilde C$.
\end{cor}
\begin{proof}
Since $\lambda^K \le N^{-1}$, the right side of~\eqref{eq:main} becomes
\beqn
\begin{split}
C_* \! \left(\frac {K+1}{\sqrt{N}} + \frac{\lambda^{K+1}}{1-\lambda}\right)  + \widetilde C\sqrt N \lambda^K
& \le
C_* \! \left(\frac{\log N + 1}{|{\log\lambda}|} + \frac{\lambda}{\sqrt N(1-\lambda)}\right)\!\frac 1{\sqrt{N}}  + \widetilde C\frac 1{\sqrt{N}}.
\end{split}
\eeqn
This yields the claim, as $\log N > 1$.
\end{proof}
As a matter of fact, the applications in Section~\ref{sec:applications} satisfy the conditions of Corollary~\ref{cor:exponential}. As a trivial example, for the map $T(x) = 2x\pmod 1$ of the interval~$[0,1]$ and a Lipschitz continuous function $f:[0,1]\to\bR$, we obtain $\lambda = \frac12$; see~Section~\ref{sec:2xmod1}. The convergence rate $\log N/\sqrt N$ is optimal up to the logarithmic correction in the numerator, which is due to the numerator of the term $K/\sqrt{N}$ in \eqref{eq:main}.  

We emphasize that the constants in Theorem~\ref{thm:main} and Corollary~\ref{cor:exponential} have explicit expressions in terms of the assumptions. In particular, the bounds do not depend on the covariance~$\Sigma$ at all. In the examples of Section~\ref{sec:applications} we also see that the bound scales as $d^3$ with increasing dimension, due to the dominating contribution of $C_*$.

\subsection{Discrete time: univariate case}
While Theorem~\ref{thm:main} readily applies in the univariate case $d=1$, we also have the following complementary result, where the assumed smoothness of the test function~$h$ is relaxed to Lipschitz continuity. The expense is that (A2) is replaced with an assumption involving an entire class of less smooth functions. 

Thus, let $Z \sim \cN(0,\sigma^2)$ be a random variable with normal distribution of mean zero and variance~$\sigma^2$. The Wasserstein distance of $W$ and $Z$ is defined as
\beqn
d_\mathscr{W}(W,Z) = \sup_{h\in\mathscr{W}}|\mu(h(W)) - \Phi_{\sigma^2}(h)|,
\eeqn
where
\beqn
\mathscr{W} = \{h:\bR\to\bR\,:\,|h(x) - h(y)| \le |x-y|\}
\eeqn
is the class of all $1$-Lipschitz functions. Theorem~\ref{thm:main-1d} below provides a bound on $d_\mathscr{W}(W,Z)$. 

Note that replacing $\mathscr{W}$ with the class of step functions
$
\mathscr{K} = \{1_{(-\infty,x]}\,:\,x\in\bR\}
$
result in the Kolmogorov distance $d_\mathscr{K}(W,Z)$, which is known to satisfy
\beqn
d_\mathscr{K}(W,Z) \le (2\pi^{-1})^{\frac14} \sigma^{-\frac12}\, d_\mathscr{W}(W,Z)^{\frac12}.
\eeqn
Hence, Theorem~\ref{thm:main-1d} also implies a Berry--Esseen type bound on the Kolmogorov distance, albeit one far from optimal.\footnote{To give an example, if $S_N$ is the displacement of a simple symmetric random walk at time~$N$, then the rate of convergence of $\frac 1{\sqrt{N}}S_N$ to the standard normal distribution is $\frac 1{\sqrt{N}}$ for both metrics.} On the other hand, there is no ``uniform'' bound in the opposite direction, so starting from a bound on $d_\mathscr{K}$ would not yield a bound on $d_\mathscr{W}$.

Moreover, the authors of~\cite{BarbourHall_1985} make the following argument in favor of metrics more regular than the Kolmogorov distance: ``Smooth metrics [\dots]\@ would seem more natural in a theoretical setting than the unsmooth uniform metric, whose principal importance is in applications to the traditional theory of significance testing. And even here, the two-point action space (accept or reject the null hypothesis), which
makes estimates [of $\mu\{W\le x\}$] important in the construction of tests, leads not only to mathematical awkwardness, but also to an analogous philosophical awkwardness, since a very small difference in observation may make a big difference in the action taken.''

We also point out that the Wasserstein distance has an interesting connection with the concept of coupling: if~$X$ and~$Y$ are random variables, then
$
d_\mathscr{W}(X,Y) = \inf_{(\hat X,\hat Y)}\bE|\hat X-\hat Y|,
$
where the infimum is taken over all random pairs $(\hat X,\hat Y)$ on a single probability space such that $\hat X$ has the distribution of $X$ and $\hat Y$ has the distribution of $Y$.

Before stating the theorem, let us recall that
a function $g:\bR\to\bR$ is said to be absolutely continuous if it has a derivative~$g'$ almost everywhere, the derivative is locally integrable, and
\beqn
g(y) = g(x) + \int_x^y g'(t)\,dt
\eeqn
for all real numbers $x\le y$.

\begin{thm}\label{thm:main-1d} Let \(f: \, X \to \bR\) be a bounded measurable function with \(\mu(f) = 0\). Fix integers $N>0$ and $0\le K < N$. Suppose that the following conditions are satisfied:

\begin{itemize}
\item[\bf (B1)] There exist constants $C_2 > 0$ and $C_4 > 0$, and a non-increasing function \(\rho : \, \bN_0 \to \bR_+\) with \(\rho(0) = 1\) and \( \sum_{i=1}^{\infty} i\rho(i) < \infty\), such that
\begin{align*}
|\mu(f\, f^k)| & \le C_2\,\rho(k)
\\
|\mu(f\, f^l f^m f^n)| & \le C_4 \min\{\rho(l),\rho(n-m)\}
\\
|\mu(f\, f^l f^m f^n) - \mu(f\, f^l)\mu(f^m f^n)| & \le C_4\,\rho(m-l)
\end{align*}
hold whenever $k\ge 0$ and $0\le l\le m \le n<N$.

\smallskip
\item[\bf (B2)] There exists a function \(\tilde\rho : \, \bN_0 \to \bR_+\) such that, given a differentiable $A:\bR\to\bR$ with $A'$ absolutely continuous and $\max_{0\le k\le 2}\|A^{(k)}\|_\infty \le 1$,
\begin{align*}
|\mu( f^n A(W^n) ) | \le \tilde\rho(K)
\end{align*}
holds for all $0\le n < N$.

\smallskip
\item[\bf (B3)] $f$ is not a coboundary. 

\end{itemize}
Then
\beq\label{eq:sigma^2}
\sigma^2 = \mu(f\,f) + 2\sum_{n=1}^{\infty} \mu(f\, f^n)
\eeq
is strictly positive and finite.
Moreover, if $Z \sim \cN(0,\sigma^2)$ is a random variable with normal distribution of mean zero and variance~$\sigma^2$, then
\begin{align}\label{eq:main-1d}
d_\mathscr{W}(W,Z) \le C_{\#} \! \left(\frac{K+1}{\sqrt{N}} + \sum_{i= K+1}^{\infty}  \rho(i) \right) + C_{\#}'\sqrt N\tilde\rho(K),
\end{align}
where
\beqn
C_{\#} = 11 \max\{\sigma^{-1},\sigma^{-2}\} \max\{C_2,\sqrt{C_4}\} (1 + \Vert f \Vert_{\infty}) \sum_{i = 0}^{\infty}(i+1) \rho(i)
\eeqn
and
\beqn
C_{\#}' = 2\max\{1,\sigma^{-2}\}
\eeqn
are independent of $N$ and $K$. 
\end{thm}

Also this theorem is proved in Sections~\ref{sec:preliminary} and~\ref{sec:proofs}. Similar comments apply to Assumptions (B1)--(B3) as to (A1)--(A3). In practice, Theorem~\ref{thm:main-1d} is applied in the same fashion as Theorem~\ref{thm:main}. For instance, a corollary analogous to Corollary~\ref{cor:exponential} is obtained verbatim. 

Let us reiterate that the qualitative difference between~(A2) and~(B2) is dictated by passing to less regular test functions. (The Wasserstein distance is defined in terms of $1$-Lipschitz test functions, not smooth ones.) Note that the absolute continuity of $A'$ implies the existence of a locally integrable $A''$, which is also assumed to be bounded in~(B2).

There is, however, a notable difference between the conclusions of the theorems: the bound in~\eqref{eq:main} of Theorem~\ref{thm:main} is independent of the covariance~$\Sigma$, while the bound in~\eqref{eq:main-1d} of Theorem~\ref{thm:main-1d} depends explicitly on the variance~$\sigma^2$. This, too, is due to the choice of metric; as we see from Theorem~\ref{thm:main} in the univariate case $d=1$, the smooth metric is insensitive to the size of the limit variance. It is not possible to entirely rid of $\sigma^2$ in Theorem~\ref{thm:main-1d}: even in the case of independent and identically distributed random variables, the bound on the Wasserstein distance depends on how the variance compares with higher absolute moments. Another example of the same phenomenon is given by the classical Berry--Esseen theorem for i.i.d.\@ variables and the Kolmogorov distance. 

Continuing the ongoing discussion, we point out that also the multivariate result, Theorem~\ref{thm:main}, could have been formulated in terms of slightly weaker regularity assumptions on the test function~$h$. Namely, from the recent proof of \cite[Proposition 2.1]{Gaunt_2016} it can be seen that sufficient conditions on~$h$ for our proofs are the following: $h\in C^2(\bR^d,\bR)$, $D^2 h$ is Lipschitz (or just uniformly) continuous, and $\max_{0\le k\le 2}\|D^k h\|_\infty<\infty$. But this would amount to two changes: (1) in place of (A2) we would need $|\mu(f^n\cdot \nabla A(W^n))|\le \tilde\rho(K)$ for all $A\in C^3(\bR^d,\bR)$ with $\|A\|_{C^3}\le 1$, analogously to (B2); and (2) the upper bound in~\eqref{eq:main} would gain an additional factor $O(\Sigma^{-1/2})$. Again we see the same phenomenon that the bound becomes sensitive to the limit covariance when the regularity of the metric decreases. Moreover, note that~$h$ is not required to be bounded in Theorem~\ref{thm:main}. Because the benefit of relaxing the regularity of~$h$ seems small compared to the cost, we have chosen to assume three bounded derivatives.

Theorems~\ref{thm:main} and~\ref{thm:main-1d} are derived for abstract dynamical systems by implementing Stein's method. We have made an effort to formulate both results in a way that would make them as easily applicable in practice as possible. In Section~\ref{sec:preliminary} we prove preliminary theorems, which imply the main results above, as shown in Section~\ref{sec:proofs}.
In the proofs the reader will find more detailed bounds on the error of normal approximation, in both the multivariate and univariate case. However, we have chosen to keep Theorems~\ref{thm:main} and~\ref{thm:main-1d} as simple as possible.

\subsection{Continuous time}
Consider a semiflow $\psi^t:X\to X$, $t\ge 0$, preserving a probability measure~$\mu$ on~$X$.
\beqn
f^t = f\circ\psi^t
\eeqn
Given $T\ge 1$, we write 
\beqn
V = V(T) = \frac{1}{\sqrt{T}} \int_0^T f^s\,ds.
\eeqn
Using a standard idea, one can immediately deduce versions of Theorems~\ref{thm:main} and~\ref{thm:main-1d} for~$V(T)$ by considering the time-one map $\cT = \psi^1$ and the observable $F:X\to\bR^d$ given by
\beqn
F = \int_0^1 f^s\,ds,
\eeqn
as follows: Denoting $N = \lfloor T\rfloor$, we have
\beqn
\begin{split}
V(T) & = \frac{\sqrt N}{\sqrt T}\frac1{\sqrt N}\int_0^N f^s\,ds + \frac{1}{\sqrt T}\int_N^T f^s\,ds
\\
& = \sqrt{1 - \frac{T-N}{T}}\frac1{\sqrt N}\sum_{n=0}^{N-1}\int_0^1 f\circ \psi^{n+s}\,ds + \frac{1}{\sqrt T}\int_N^T f^s\,ds
\\
& = a(T)\frac1{\sqrt N}\sum_{n=0}^{N-1}F\circ \cT^n + \frac{1}{\sqrt T}\int_N^T f^s\,ds,
\end{split}
\eeqn
where
\beqn
a = a(T) = \sqrt{1 - \frac{T-N}{T}}
\eeqn
satisfies
\beqn
|a(T) -1| \le  \frac{T-N}{T}.
\eeqn
We can estimate the rate of convergence of $V(T)$ to $\cN(0,\Sigma)$ in the limit $T\to\infty$, if we can do the same for 
\beqn
W = W(N) = \frac1{\sqrt N}\sum_{n=0}^{N-1}F^n
\quad\text{where} \quad F^n = F\circ \cT^n,
\eeqn
in the limit $N\to\infty$. Namely,
\beqn
\begin{split}
|(V - W)_\alpha| & \le |a-1|\left|\frac1{\sqrt N}\sum_{n=0}^{N-1}F_\alpha\circ \cT^n\right| +  \left|\frac{1}{\sqrt T}\int_N^T f^s_\alpha\,ds\right| 
\\
& \le \frac{T-N}{T}\sqrt N \|f_\alpha\|_\infty + \frac{T-N}{\sqrt T}\|f_\alpha\|_\infty \le \frac{2\|f\|_\infty}{\sqrt T} 
\end{split}
\eeqn
for all $\alpha\in\{1,\dots,d\}$. It follows that
\beqn
|h(V) - h(W)| \le \frac{2d \|\nabla h\|_\infty\|f\|_\infty}{\sqrt T}
\eeqn
holds for any Lipschitz continuous function $h:\bR^d\to\bR$, which results in the final estimate
\beq\label{eq:VW}
|h(V) - \Phi_\Sigma(h)| \le |h(W) - \Phi_\Sigma(h)| + \frac{2d \|\nabla h\|_\infty\|f\|_\infty}{\sqrt T}.
\eeq
Here Theorem~\ref{thm:main} or~\ref{thm:main-1d} can then be used to determine a bound on $|h(W) - \Phi_\Sigma(h)|$, with
\beq\label{eq:Sigma_time-one}
\Sigma = \mu(F\otimes F) + \sum_{n=1}^{\infty} (\mu(F^n \otimes F) + \mu(F \otimes F^n)).
\eeq

In spite of the above, depending on the situation, it may be preferable to work directly with the semiflow~$\psi^t$ and the observable~$f$, which is how Theorem~\ref{thm:main_flow} below is formulated. A similar result could be proved identically to Theorem~\ref{thm:main}, but to save space we deduce it directly from the latter. For that reason we introduce the quantity 
\beq\label{eq:V^t}
\begin{split}
%
V^t = V^t(T,K) = \frac{1}{\sqrt{\lfloor T\rfloor}} \int 1_{[0,\lfloor t\rfloor -K]}(s) f^{s}\,ds + \frac{1}{\sqrt{\lfloor T\rfloor}} \int 1_{[\lfloor t\rfloor+K+1,\lfloor T\rfloor]}(s) f^{s}\,ds.
\end{split}
\eeq
Here it is understood that the indicator function $1_{[a,b]}$ vanishes identically if $a>b$.

\begin{thm}\label{thm:main_flow} Let \(f: \, X \to \bR^d\) be a bounded measurable function with \(\mu(f) = 0\). Let \(h: \, \bR^d \to \bR\) be three times differentiable with \(\Vert D^k h \Vert_{\infty} < \infty\) for \(1 \le k \le 3\). Fix a real number $T\ge 1$ and an integer $0< K < T$. Suppose that the following conditions are satisfied:

\begin{itemize}
\item[{\bf (C1)}] There exist constants $C_2 > 0$ and $C_4 > 0$, and a non-increasing function \(\rho : \, \bR_+ \to \bR_+\) with \(\rho(0) = 1\) and \( \int_0^\infty t\,\rho(t)\,dt < \infty\), such that
\begin{align*}
|\mu(f_{\alpha} f^r_{\beta})| & \le C_2\,\rho(r)
\\
|\mu(f_{\alpha} f_{\beta}^s f_{\gamma}^t f_{\delta}^u)| & \le C_4 \min\{\rho(s),\rho(u-t)\}
\\
|\mu(f_{\alpha} f_{\beta}^s f_{\gamma}^t f_{\delta}^u) - \mu(f_{\alpha} f_{\beta}^s)\mu(f_{\gamma}^t f_{\delta}^u)| & \le C_4\,\rho(t-s)
\end{align*}
hold whenever $r\ge 0$; $0\le s\le t \le u \le T$; $\alpha,\beta,\gamma,\delta\in\{\alpha',\beta'\}$ and $\alpha',\beta'\in\{1,\dots,d\}$.

\smallskip
\item[{\bf (C2)}] There exists a function \(\tilde\rho : \, \bN_0 \to \bR_+\) such that
\begin{align*}
|\mu( f^t \cdot  \nabla h(w + V^t \tau) ) | \le \tilde\rho(K)
\end{align*}
$0\le t < \lfloor T\rfloor$, $0\le \tau\le 1$ and $w\in\bR^d$.

\smallskip
\item[{\bf (C3)}] Assume that the matrix\footnote{The symmetric $d\times d$ matrix is well defined and positive semidefinite.}
\beq\label{eq:Sigma_flow}
\Sigma = \int_0^\infty (\mu(f^t \otimes f) + \mu(f \otimes f^t))\,dt
\eeq
is positive definite.
\end{itemize}
Then
\beqn\label{eq:main_flow}
|\mu(h(V)) - \Phi_{\Sigma}(h)| \le 6C_*\!\left(\frac {K+1}{\sqrt{T}} + \sum_{i= K}^{\infty}  \rho(i)\right)  + \sqrt{T} \tilde\rho(K) + \frac{2d \|\nabla h\|_\infty\|f\|_\infty}{\sqrt T},
\eeqn
where $C_*$ has the same exact expression as in~\eqref{eq:main_constant}.
\end{thm}

We show in Section~\ref{sec:proofs_flow} how Theorem~\ref{thm:main_flow} is obtained from Theorem~\ref{thm:main}. Similarly, a continuous-time version of the univariate Theorem~\ref{thm:main-1d} can be obtained. We leave that to the interested reader.

Let us also make a small final remark concerning the covariance matrix $\Sigma$ and Assumption~(C3). For the time-one map, the expression of $\Sigma$ is given in~\eqref{eq:Sigma_time-one} before Theorem~\ref{thm:main_flow}. This, of course, coincides with the expression in~\eqref{eq:Sigma_flow}: elementary computations show
\beqn
\mu(F\otimes F) = \int_0^1 (1-t)\, \mu(f^t\otimes f)\,dt + \int_0^1 (1-t)\, \mu(f\otimes f^t)\,dt
\eeqn
and, for $n\ge 1$,
\beqn
\begin{split}
\mu(F^n\otimes F) & = \int_0^1 t\, \mu(f^{n-1+t}\otimes f)\,dt + \int_0^1 (1-t)\, \mu(f^{n+t}\otimes f)\,dt
\\
\mu(F\otimes F^n) & = \int_0^1 t\, \mu(f\otimes f^{n-1+t})\,dt + \int_0^1 (1-t)\, \mu(f\otimes f^{n+t})\,dt.
\end{split}
\eeqn
Summing these yields the claim.
In particular, Assumption~(C3) fails if and only if $F$ is a coboundary for the time-one map.


\section{Review of Stein's method}\label{sec:Stein}

In this section we give a brief overview of Stein's method, sufficient for understanding our paper. We refer the reader to Stein's original paper~\cite{Stein_1972} and the references \cite{BarbourChen_2005,Chen_etal_2011, Ross_2011, Chatterjee_2014} for more detailed treatments of the subject. For conceptual simplicity, we begin with the case of univariate normal approximation.

\subsection{Univariate normal approximation}
The results of our paper concern convergence to a normal distribution, the error being measured using the Wasserstein distance. Here we review Stein's method in that case. 

Below, $\cN(0,\sigma^2)$ will stand for the centered normal distribution with variance $\sigma^2>0$. Let $\Phi_{\sigma^2}(h)$ denote the expectation of a function $h:\bR\to\bR$ with respect to $\cN(0,\sigma^2)$, the centered normal distribution with variance $\sigma^2$. That is,
\beqn
\Phi_{\sigma^2}(h) = \frac{1}{\sqrt{2\pi}\sigma}\int_{-\infty}^\infty e^{-t^2/2\sigma^2} h(t)\,dt.
\eeqn

To proceed, let us recall a lemma due to Stein:

\begin{lem}\label{lem:Stein}
A random variable $W$ has the distribution $\cN(0,\sigma^2)$ if and only if
\beqn
\bE[\sigma^2 A'(W) - WA(W)] = 0
\eeqn
for all absolutely continuous functions $A:\bR\to\bR$ satisfying $\Phi_{\sigma^2}(|A'|)<\infty$.
\end{lem}

The proof of Lemma~\ref{lem:Stein} can be found, e.g., in~\cite{Chen_etal_2011}. What is important is that Stein's lemma characterizes normal distribution by a simple equation involving a functional of the random variable~$W$. This suggests that if $\bE[\sigma^2 A'(W) - WA(W)] \approx 0$ for a sufficiently large class of functions~$A$, then the distribution of $W$ should be close to~$\cN(0,\sigma^2)$. In order to make this precise, Stein introduced the equation
\beq\label{eq:Stein}
\sigma^2 A'(w) - wA(w) = h(w) - \Phi_{\sigma^2}(h)
\eeq
with the following idea: Suppose that, for each test function $h$ belonging to some class~$\mathscr{H}$, the Stein equation~\eqref{eq:Stein} has a solution $A$ belonging to another class of functions~$\mathscr{A}$. Then, for $Z\sim\cN(0,\sigma^2)$,
\beqn
d_\mathscr{H}(W,Z) \equiv \sup_{h\in\mathscr{H}} |\bE h(W) - \Phi_{\sigma^2}(h)| \le \sup_{A\in\mathscr{A}} |\bE[\sigma^2 A'(W) - WA(W)]|.
\eeqn
Thus, the distance between the distribution of $W$ and $\cN(0,\sigma^2)$ can be bounded if the right side can be bounded, and this only involves working with the distribution of $W$. Of course, different choices of~$\mathscr{H}$ yield different classes~$\mathscr{A}$. For the Wasserstein distance ($\mathscr{H} = \mathscr{W}$) the following is known:

\begin{lem}\label{lem:univ_lemma}
Let $Z$ be a random variable with distribution $\cN(0,\sigma^2)$ and $W$ any random variable. Then
\beq\label{eq:univ_bound}
d_\mathscr{W}(W,Z) \le \sup_{A\in\mathscr{F}_{\sigma^2}} |\bE[\sigma^2 A'(W) - WA(W)]|
\eeq
where $\mathscr{F}_{\sigma^2}$ is the class of all differentiable functions $A:\bR\to\bR$ with an absolutely continuous derivative, satisfying the bounds
\beqn
\|A\|_\infty \le 2, 
\quad
\|A'\|_\infty \le \sqrt{2/\pi}\, \sigma^{-1}
\quad\text{and}\quad
\|A''\|_\infty \le 2 \sigma^{-2}.
\eeqn
\end{lem}

\begin{proof}[Sketch of proof]
The proof for $\sigma^2 = 1$ (see~\cite{Chen_etal_2011}) is based on showing that the solution
\beqn
A_1(w) = e^{w^2/2}\int_w^\infty e^{-t^2/2}(\Phi_1(h_1)-h_1(t))\,dt
\eeqn
to
\beqn
A_1'(w) - w A_1(w) = h_1(w) - \Phi_1(h_1)
\eeqn
with an absolutely continuous test function $h_1$ satisfies the bounds
\beqn
\|A_1\|_\infty \le 2\|h_1'\|_\infty, 
\quad
\|A_1'\|_\infty \le \sqrt{2/\pi}\,\|h_1'\|_\infty
\quad\text{and}\quad
\|A_1''\|_\infty \le 2\|h_1'\|_\infty .
\eeqn
The general case can be reduced to this by changes of variables as follows.
Suppose $\sigma^2 > 0$ and $h$ is an absolutely continuous test function. Denoting
\beqn
h_1(w) = h(\sigma w), 
\eeqn
we have
$
\Phi_{\sigma^2}(h) = 
\Phi_1(h_1).
$
Defining $A_1$ as above, the function\footnote{In fact, $A(w) = \sigma^{-2}e^{w^2/2\sigma^2}\int_{w}^\infty e^{-t^2/2\sigma^2}(\Phi_{\sigma^2}(h)-h(t))\,dt$.}
\beqn
A(w) = \sigma^{-1}A_1(\sigma^{-1} w)
\eeqn
solves~\eqref{eq:Stein} and, by the estimates above,
\beqn
\|A\|_\infty \le 2\|h'\|_\infty, 
\quad
\|A'\|_\infty \le \sqrt{2/\pi}\,\sigma^{-1}\|h'\|_\infty
\quad\text{and}\quad
\|A''\|_\infty \le  2\sigma^{-2}\|h'\|_\infty.
\eeqn
Finally, any $h\in\mathscr{W}$ is absolutely continuous with $\|h'\|_\infty\le 1$.
\end{proof}

Thus, in the univariate case the task is to bound the right side of~\eqref{eq:univ_bound}.

\subsection{Multivariate normal approximation}
Stein's method for multivariate normal approximation is similar, so we only record the bits of the theory relevant to our work. There is one essential difference to the univariate case: the test function~$h$ is generally required to be more regular in dimensions greater than one.

Let the matrix~\(\Sigma \in \bR^{d\times d}\) be symmetric and positive definite. Denote by~$\phi_\Sigma$ and~$\Phi_\Sigma$ the density and distribution function of the $d$-dimensional normal distribution with mean $0$ and covariance $\Sigma$, respectively. Given a test function \(h:\bR^d\to\bR\), define
\begin{align}\label{eq:A}
&A(w)= -\int_0^{\infty} \! \left\{\int_{\bR^d} h(e^{-s}w + \sqrt{1 - e^{-2s}}\,z)\,\phi_\Sigma(z)\,dz - \Phi_{\Sigma}(h)\right\} \, ds,
\end{align}
where $\Phi_\Sigma(h)$ stands for the expectation of $h$ with respect to the above normal distribution.
Then, we have the following result; see~\cite{Barbour_1990,Gotze_1991,GoldsteinRinott_1996,Gaunt_2016}.

\begin{lem}\label{lem:steinmv} Let \(h: \, \bR^d \to \bR\) be three times differentiable with \(\Vert D^k h \Vert_{\infty} < \infty\) for \(1 \le k \le 3\). Then,  \(A \in C^3(\bR^d,\bR)\), and \(A\) solves the Stein equation
\begin{align}\label{eq:steinmv}
\tr  \Sigma D^2A(w) - w \cdot  \nabla A(w) = h(w) - \Phi_{\Sigma}(h).
\end{align}
Moreover, the partial derivatives of \(A\) satisfy the bounds 
 \begin{align*}
\| \partial_1^{t_1}\cdots \partial_d^{t_d} A \|_\infty \le k^{-1} \| \partial_1^{t_1}\cdots \partial_d^{t_d} h \|_\infty
\end{align*}
whenever $t_1+\cdots+t_d = k$, $1\le k\le 3$.
\end{lem}

Thus, in the multivariate case the task is to bound
\beqn
|\bE h(W) - \Phi_{\Sigma}(h)| = |\bE[\tr  \Sigma D^2A(W) - W \cdot  \nabla A(W)]|.
\eeqn
 
Note the formal difference to the univariate case that in the present, multivariate, case the Stein equation is of second order. Note also that the bounds on $A$ in Lemma~\ref{lem:steinmv} are independent of the covariance~$\Sigma$, whereas in Lemma~\ref{lem:univ_lemma} they depend on the variance~$\sigma^2$. This has to do with the additional regularity required of~$h$ in the multivariate case; see the earlier discussion after Theorem~\ref{thm:main-1d}, in particular the comments pertaining to~\cite{Gaunt_2016}.

\subsection{Other target distributions}
In this section we briefly underline that Stein's method is not limited to normal approximation. It would be interesting to find applications of some of the ideas mentioned here to the theory of dynamical systems.

The general outline of the method remains the same for other target distributions, and is as follows: Let $W$ be any random variable and $Z$ a random variable with some specific distribution~$\nu_Z$. The task is compare the distributions of $W$ and $Z$ when integrated against test functions~$h$ in a suitable class of functions~$\mathscr{H}$, that is, to bound
\beqn
\sup_{h\in\mathscr{H}}|\bE h(W)- \nu_Z(h)|.
\eeqn  

The strategy to accomplish this is to determine a Stein operator $\cS$ which characterizes the distribution of~$Z$ in the sense that $\bE[\cS A(W)]=0$ holds for all $A$ in some class of functions $\mathscr{A}$ if and only if $W \stackrel d= Z$. Furthermore, the classes $\mathscr{H}$ and $\mathscr{A}$ should be compatible in the sense that, given $h\in\mathscr{H}$, the Stein equation
\beqn
\cS A(w) = h(w) - \nu_Z(h)
\eeqn
has a solution $A\in\mathscr{A}$. Since the latter satisfies 
\beqn
\bE h(W) - \nu_Z(h) = \bE[\cS A(W)],
\eeqn
the original task is reduced to bounding 
\beqn
\sup_{A\in\mathscr{A}}|\bE[\cS A(W)]|,
\eeqn
which only involves the random variable~$W$. For completeness, we list the Stein operators of some common target distributions:
\begin{itemize} 
\item 
Poisson distribution with mean $\lambda>0$ \cite{Chen_1975,ArratiaGoldsteinGordon_1990}:
\beqn
\cS A(w) = \lambda A(w+1)- wA(w).
\eeqn
This operator is exactly what the dynamical systems works \cite{Denker_etal_2004, Psiloyenis_2008, GordinDenker_2012, Haydn_2013, HaydnYang_2016}, mentioned in the introduction, build on.

\medskip
\item
Exponential distribution with mean one \cite{FulmanRoss_2013}:
\beqn
\cS A(w) = wA'(w)-(w-1)A(w).
\eeqn

\medskip
\item
Binomial distribution of~$n$ independent experiments with success probability~$p$ \cite{Ehm_1991}:
\beqn
\cS A(w) = p(n - w)A(w + 1) - (1-p)wA(w).
\eeqn

\medskip
\item
Gamma distribution $\Gamma(r,\lambda)$\cite{Luk_1994}:
\beqn
\cS A(w) = wA''(w)+(r-\lambda w)f'(w).
\eeqn
 

\end{itemize}

\noindent Applications of Stein's method involving these and other target distributions can be readily found in the vast probability theory literature. Stein's method has also applications at least to functional limit theorems~\cite{Barbour_1990}, lower bounds on the error of normal approximation~\cite{HallBarbour_1984}, large deviation bounds~\cite{Raic_2007}, concentration inequalities~\cite{Chatterjee_2005,Chatterjee_2007,ChatterjeeDey_2010}, and transport inequalities~\cite{LedouxNourdinPeccati_2015}.


\section{Preliminary results in discrete time}\label{sec:preliminary}

In this section we prove preliminary versions of Theorems~\ref{thm:main} and~\ref{thm:main-1d}, obtaining explicit upper bounds in terms of the assumptions.  

\subsection{Statements}

For the convenience of the reader, let us recall from the beginning of Section~\ref{sec:results} that we denote $W = \frac{1}{\sqrt{N}}\sum_{k=0}^{N-1} f^k$ and \(W^n = W - \frac{1}{\sqrt{N}} \sum_{k\in [n]_K} f^k\), where furthermore \([n]_K =  \{ k \in \bN_0 \cap [0,N-1] \, : \, |k-n| \le K \}\).

The following result is a preliminary version of~Theorem~\ref{thm:main}. Again for the reader's convenience, we repeat some of the assumptions:

\begin{thm}\label{thm:pre} Let \(f: \, X \to \bR^d\) be a bounded measurable function with \(\mu(f) = 0\). Let $A\in C^3(\bR^d,\bR)$ be a given function satisfying \(\| D^k A \|_{\infty} < \infty\) for \(1 \le k \le 3\). Fix integers $N>0$ and $0\le K< N$.  Suppose that the following conditions are satisfied:

\begin{itemize}
\item[(A1')] There exist constants $C_2 > 0$ and $C_4 > 0$, and a non-increasing function \(\rho : \, \bN_0 \to \bR_+\) with \(\rho(0) = 1\) and \( \sum_{i=1}^{\infty} \rho(i) < \infty\) \footnote{Assumption (A1') differs from (A1) only in that the latter also requires $\sum_{i=1}^\infty i\rho(i)<\infty$.}, such that
\begin{align*}
|\mu(f_{\alpha} f^k_{\beta})| & \le C_2\,\rho(k)
\\
|\mu(f_{\alpha} f_{\beta}^l f_{\gamma}^m f_{\delta}^n)| & \le C_4 \min\{\rho(l),\rho(n-m)\}
\\
|\mu(f_{\alpha} f_{\beta}^l f_{\gamma}^m f_{\delta}^n) - \mu(f_{\alpha} f_{\beta}^l)\mu(f_{\gamma}^m f_{\delta}^n)| & \le C_4\,\rho(m-l)
\end{align*}
hold whenever $k\ge 0$; $0\le l\le m \le n< N$; $\alpha,\beta,\gamma,\delta\in\{\alpha',\beta'\}$ and $\alpha',\beta'\in\{1,\dots,d\}$.

\smallskip
\item[(A2')] There exists a function \(\eta: \, \bN_0^2 \to \bR_+\) such that
\begin{align*}
\left|\mu\!\left(\frac{1}{\sqrt{N}}\sum_{n=0}^{N-1} f^n \cdot  \nabla A(W^n)\right)\right|  \leq \eta(N,K).
\end{align*}

\smallskip
\item[(A3)] $f$ is not a coboundary in any direction. 

\end{itemize}
Then $\Sigma$ in~\eqref{eq:Sigma} is a well-defined, symmetric, positive-definite, $d\times d$ matrix; and
\begin{align*}
&|\mu(\tr \Sigma D^2A(W) - W \cdot  \nabla A(W))| \\
&\qquad \le d^3 C_2\Vert f \Vert_{\infty} \Vert D^3 A \Vert_{\infty}\frac{2K+1}{\sqrt{N}} \left(\rho(0) + 2\sum_{i = 1}^{2K} \rho(i) \right)
 \\
&\qquad\quad+ 2 d^2C_2\Vert D^2 A \Vert_{\infty}  \left( \sum_{i= K+1}^{\infty}  \rho(i)  + \frac1N\sum_{i= 1}^{K} i \rho(i) \right)\\
&\qquad\quad+ 11d^2\max\{C_2,\sqrt{C_4}\} \Vert D^2 A\Vert_{\infty}\frac{\sqrt{K+1}}{\sqrt{N}}\sqrt{\sum_{i=0}^{N-1} (i+1)\rho(i)}\, \\
&\qquad\quad+ \eta(N,K).
\end{align*}
\end{thm}

The following result is a preliminary version of Theorem~\ref{thm:main-1d}. It will be proved in Section~\ref{sec:proof_pre-1d}.

\begin{thm}\label{thm:pre-1d}
Let \(f: \, X \to \bR\) be a bounded measurable function with \(\mu(f) = 0\). Let $A\in C^1(\bR,\bR)$ be a given function with absolutely continuous~$A'$, satisfying $\|A^{(k)}\|_\infty < \infty$ for $0\le k\le 2$. Fix integers $N>0$ and $0\le K<N$. Suppose that the following conditions are satisfied:
\begin{itemize}
\item[(B1')] There exist constants $C_2 > 0$ and $C_4 > 0$, and a non-increasing function \(\rho : \, \bN_0 \to \bR_+\) with \(\rho(0) = 1\) and \( \sum_{i=1}^{\infty} \rho(i) < \infty\) \footnote{Assumption (B1') differs from (B1) only in that the latter also requires $\sum_{i=1}^\infty i\rho(i)<\infty$.}, such that
\begin{align*}
|\mu(f\, f^k)| & \le C_2\,\rho(k)
\\
|\mu(f\, f^l f^m f^n)| & \le C_4 \min\{\rho(l),\rho(n-m)\}
\\
|\mu(f\, f^l f^m f^n) - \mu(f\, f^l)\mu(f^m f^n)| & \le C_4\,\rho(m-l)
\end{align*}
hold whenever $k\ge 0$ and $0\le l\le m \le n<N$.

\smallskip
\item[(B2')] There exists a function \(\eta: \, \bN_0^2 \to \bR_+\) such that
\begin{align*}
\left|\mu\!\left(\frac{1}{\sqrt{N}}\sum_{n=0}^{N-1} f^nA(W^n)\right)\right| \le \eta(N,K).
\end{align*}

\smallskip
\item[(B3)] $f$ is not a coboundary. 

\end{itemize}

Then~$\sigma^2$ in~\eqref{eq:sigma^2} is strictly positive and finite; and
\begin{align*}
&|\mu(\sigma^2 A'(W)-WA(W))| \\
&\qquad \le \frac12 C_2\Vert f \Vert_{\infty} \Vert A'' \Vert_{\infty}\frac{2K+1}{\sqrt{N}} \left(\rho(0) + 2\sum_{i = 1}^{2K} \rho(i) \right)
 \\
&\qquad\quad+ 2C_2\Vert A' \Vert_{\infty}  \left( \sum_{i= K+1}^{\infty}  \rho(i)  + \frac1N\sum_{i= 1}^{K} i \rho(i) \right)\\
&\qquad\quad+ 11\max\{C_2,\sqrt{C_4}\} \Vert A'\Vert_{\infty}\frac{\sqrt{K+1}}{\sqrt{N}}\sqrt{\sum_{i=0}^{N-1} (i+1)\rho(i)}\, \\
&\qquad\quad+ \eta(N,K).
\end{align*}

\end{thm}
The proof of Theorem~\ref{thm:pre-1d} is nearly identical to that of Theorem~\ref{thm:pre}. The minor differences are that in $d=1$ the Stein equation can be reduced to a first order differential equation (hence the lower order derivatives of~$A$ in the upper bound) and Taylor's theorem can be used for absolutely continuous functions instead of continuously differentiable ones. These will be detailed in Section~\ref{sec:proof_pre-1d}.

\subsection{Proof of Theorem~\ref{thm:pre}}\label{sec:proof_pre}
By Assumption (A1'), the series
\beqn
\Sigma = \mu(f \otimes f) + \sum_{n=1}^{\infty} (\mu(f^n \otimes f) + \mu(f \otimes f^n))
\eeqn
is absolutely convergent componentwise. Let us remark that the slightly stronger Assumption~(A1) of Theorem~\ref{thm:main} also guarantees $\lim_{N\to\infty} \mu(W\otimes W) = \Sigma$, but this is not necessary for Thoerem~\ref{thm:pre} at hand. That Assumption~(A3) implies the positive-definiteness of~$\Sigma$ has been discussed in Section~\ref{sec:results}.

Recall that the task is to obtain an upper bound on
$
|\mu(\tr \Sigma D^2A(W) - W \cdot  \nabla A(W))|,
$
where \(A \in C^3(\bR^d,\bR)\) with \( \Vert D^k A \Vert_{\infty} < \infty\), \(k=1,2,3\).
To this end, let us write
\begin{align*}
& \mu(W\cdot  \nabla A(W)) 
\\
&\qquad = \mu \! \left(\frac{1}{\sqrt{N}}\sum_{n=0}^{N-1} f^n\cdot ( \nabla A(W) -  \nabla A(W^n))\right) + \mu \! \left(\frac{1}{\sqrt{N}}\sum_{n=0}^{N-1} f^n \cdot  \nabla A(W^n)\right),
\\
&\qquad = \mu\!\left( \frac{1}{\sqrt{N}} \sum_{n=0}^{N-1} f^n \cdot (  \nabla A(W) -  \nabla A(W^n) - D^2A(W)(W-W^n)) \right) \\
& \quad\qquad + \mu \! \left( \frac{1}{\sqrt{N}} \sum_{n=0}^{N-1} f^n \cdot D^2A(W)(W-W^n) \right) 
\\
& \quad\qquad + \mu \! \left(\frac{1}{\sqrt{N}}\sum_{n=0}^{N-1} f^n \cdot  \nabla A(W^n)\right).
\end{align*}
For all vectors $u,v\in\bR^d$ and matrices $M\in\bR^{d\times d}$, the identity
\beqn
\tr(u\otimes v \, M) = \sum_{\alpha=1}^d (u\otimes v\, M)_{\alpha\alpha}  =  \sum_{\alpha,\beta=1}^d u_\alpha v_\beta M_{\beta\alpha}
 = u\cdot M^t v 
\eeqn
holds. Since $D^2 A(W)$ is symmetric and 
\beqn
W-W^n = \frac{1}{\sqrt{N}} \sum_{m\in [n]_K} f^m,
\eeqn
we thus have
\begin{align*}
 \frac{1}{\sqrt{N}} \sum_{n=0}^{N-1} f^n \cdot D^2A(W)(W-W^n) = \frac{1}{N} \sum_{n=0}^{N-1}\sum_{m \in [n]_K} \tr (f^n \otimes f^mD^2A(W)),
\end{align*}
so that
\begin{align*}
&\tr \Sigma D^2A(W) -   \frac{1}{\sqrt{N}} \sum_{n=0}^{N-1} f^n \cdot D^2A(W)(W-W^n) \\
= & \tr \Sigma D^2A(W) - \frac{1}{N} \sum_{n=0}^{N-1}\sum_{m \in [n]_K} \tr(f^n \otimes f^mD^2A(W)) \\
= & \tr \! \left(\left(  \Sigma - \frac{1}{N}\sum_{n=0}^{N-1}\sum_{m \in [n]_K} f^n \otimes f^m \right)\! D^2A(W) \right).
\end{align*}
Consequently, we arrive at the intermediate bound
\begin{align}
&|\mu(\tr \Sigma D^2A(W) - W \cdot  \nabla A(W))|\nonumber
\\
&\qquad \le \left| \mu\!\left( \frac{1}{\sqrt{N}} \sum_{n=0}^{N-1} f^n \cdot (  \nabla A(W) -  \nabla A(W^n) - D^2A(W)(W-W^n)) \right) \right| \label{eq1}\\
&\qquad\quad + \left| \mu\!\left( \tr \left(\left(  \Sigma - \frac{1}{N}\sum_{n=0}^{N-1}\sum_{m \in [n]_K} f^n \otimes f^m \right)\! D^2A(W) \right)\right) \right|\label{eq3}\\
&\qquad\quad +\left| \mu\!\left(\frac{1}{\sqrt{N}}\sum_{n=0}^{N-1} f^n \cdot  \nabla A(W^n)\right) \right| \label{eq2}.
\end{align}
By Assumption~(A2'), \eqref{eq2} is bounded by $\eta(N,K)$. It remains to obtain bounds on \eqref{eq1} and \eqref{eq3}, which is done in Sections~\ref{sec:bound1} and~\ref{sec:bound3}, respectively.

\subsubsection{Bound on \eqref{eq1}}\label{sec:bound1}
\begin{prop}\label{prop1}
The expression in~\eqref{eq1} is bounded by
\begin{align*}
d^3 C_2\Vert f \Vert_{\infty} \Vert D^3 A \Vert_{\infty}\frac{2K+1}{\sqrt{N}} \left(\rho(0) + 2\sum_{i = 1}^{2K} \rho(i)  \right).
\end{align*}
 \end{prop}
\begin{proof} Since~$A\in C^3(\bR^d,\bR)$, Taylor expansion of $\nabla A(W^n)$ at $W$ yields
\begin{align*}
&  \nabla A(W) -  \nabla A(W^n) - D^2A(W)(W-W^n)
\\
=\ & - \sum_{\beta, \gamma = 1}^d R_{\beta\gamma}(W^n)(W^n-W)_{\beta}(W^n-W)_{\gamma},
\end{align*}
where
\begin{align*}
R_{\beta\gamma} = (1 + 1_{\beta \neq \gamma}) \int_0^1 (1-u)\partial^2_{\beta\gamma} \nabla A(W + u(W^n-W)) \, du.
\end{align*}
It follows that \eqref{eq1} can be written in the form
\begin{align*}
&\left| \mu\!\left( \frac{1}{\sqrt{N}} \sum_{n=0}^{N-1} \sum_{\beta,\gamma = 1}^d f^n\cdot R_{\beta\gamma}\,(W^n-W)_{\beta}(W^n-W)_{\gamma} \right) \right|.
\end{align*}
Since
\begin{align*}
|f^n\cdot R_{\beta\gamma}| \le \sum_{1\le\alpha\le d}   |f^n_\alpha||(R_{\beta\gamma})_\alpha| \le d\|f\|_\infty \max_{1\le\alpha\le d} |(R_{\beta\gamma})_\alpha| \le d\|f\|_\infty \max_{1\le\alpha\le d} \|\partial_{\alpha\beta\gamma}^3 A\|_\infty,
\end{align*}
this yields the following bound on \eqref{eq1}:
\begin{align*}
&\le \frac{d\,\Vert f \Vert_{\infty} \Vert D^3A \Vert_{\infty}}{\sqrt{N}} \sum_{n=0}^{N-1} \sum_{\beta,\gamma=1}^d  \mu(|(W^n-W)_{\beta}(W^n-W)_{\gamma}|)
\\
&\le \frac{d\,\Vert f \Vert_{\infty} \Vert D^3A \Vert_{\infty}}{\sqrt{N}} \sum_{n=0}^{N-1} \sum_{\beta,\gamma=1}^d  \mu((W^n-W)_{\beta}^2)^\frac12\mu((W^n-W)_{\gamma}^2))^\frac12
\\
&\le \frac{d\, \Vert f \Vert_{\infty} \Vert D^3A \Vert_{\infty}}{\sqrt{N}} \sum_{n=0}^{N-1} d^2\max_{1\le \beta\le d} \mu((W^n-W)_{\beta}^2)
\\
& \le \frac{d^3\Vert f \Vert_{\infty} \Vert D^3A \Vert_{\infty}}{\sqrt N} \sum_{n=0}^{N-1} \max_{1\le \beta\le d} \mu \! \left(\left(\frac1{\sqrt N}\sum_{k\in[n]_K}f_{\beta}^k\right)^2\right)
\\
& = \frac{d^3\Vert f \Vert_{\infty} \Vert D^3A \Vert_{\infty}}{N^{3/2}} \sum_{n=0}^{N-1} \max_{1\le\beta\le d} \sum_{k,m\in [n]_K} \mu (f_{\beta}^kf_{\beta}^m).
\end{align*}
Above, the second line uses H\"older's inequality. Invoking invariance and Assumption~(A1'),
\begin{align*}
 \sum_{k,m\in [n]_K} \mu (f_{\beta}^kf_{\beta}^m)
& =  \sum_{0\le k,m< |[n]_K|} \mu (f_{\beta}^kf_{\beta}^m)
\\
& = |[n]_K| \mu(f_{\beta}f_\beta) + 2\sum_{0\le k<m<|[n]_K|} \mu (f_{\beta}^kf_{\beta}^m) 
\\
& = |[n]_K| \mu(f_{\beta}f_\beta) + 2\sum_{1\le i<|[n]_K|} (|[n]_K| - i)\mu (f_{\beta}f_{\beta}^i) 
\\
& \le  |[n]_K| C_2   \!\left( \rho(0) + 2\sum_{1\le i<|[n]_K|} \rho(i) \right)
\\
& \le (2K+1) C_2\!\left(\rho(0) + 2\sum_{1\le i\le 2K} \rho(i)\right)
\end{align*}
for all $\beta$. The last line uses $|[n]_K|\le 2K+1$. This finishes the proof of Proposition~\ref{prop1}.
\end{proof}

\subsubsection{Bound on~\eqref{eq3}}\label{sec:bound3}
\begin{prop}\label{prop3}
The expression in~\eqref{eq3} is bounded by
\beqn
\begin{split}
& 2 d^2C_2 \Vert D^2 A \Vert_{\infty}  \left( \sum_{i = K+1}^{\infty}  \rho(i)  + \frac1N \sum_{i = 1}^{K} i \rho(i) \right)
\\
+\ &
11d^2\max\{C_2,\sqrt{C_4}\} \| D^2 A\|_{\infty}\frac{\sqrt{K+1}}{\sqrt{N}}\sqrt{\sum_{i=0}^{N-1} (i+1)\rho(i)}.
\end{split}
\eeqn
\end{prop}
The rest of this section comprises the proof of Proposition~\ref{prop3}. Define
\begin{align*}
\widetilde\Sigma = \frac{1}{N}\sum_{n=0}^{N-1}\sum_{m \in [n]_K} \mu(f^n \otimes f^m),
\end{align*}
so that \eqref{eq3} has the upper bound
\begin{align}
& \mu\!\left(\left| \tr \! \left(\left(  \widetilde\Sigma - \frac{1}{N}\sum_{n=0}^{N-1}\sum_{m \in [n]_K} f^n \otimes f^m \right)\! D^2A(W) \right)\right|\right) \label{eq4} \\
& + \| { \tr (( \Sigma - \widetilde\Sigma ) D^2A ) } \|_\infty. \label{eq5}
\end{align}

\begin{lem}\label{lem5}
The expression in~\eqref{eq5} satisfies the bound
\begin{align*}
 \| { \tr (( \Sigma - \widetilde\Sigma ) D^2A ) }  \|_\infty  \le 2 d^2C_2 \Vert D^2 A \Vert_{\infty}  \left( \sum_{i= K+1}^{\infty}  \rho(i)  + \frac1N\sum_{i= 1}^{K} i \rho(i) \right).
\end{align*} 
\end{lem}

\begin{proof} 
Using invariance, an elementary computation yields
\begin{align*}
\widetilde\Sigma 
& = \mu(f \otimes f) + \frac{1}{N} \sum_{n=1}^{N-1} \sum_{\substack{m \in [n]_K \\ m < n }} \mu(f^{n-m} \otimes f) + \frac{1}{N} \sum_{n=0}^{N-2} \sum_{\substack{m \in [n]_K \\ m > n }} \mu(f \otimes f^{m-n}).
\end{align*}
Here
\begin{align*}
\sum_{n=1}^{N-1} \sum_{\substack{m \in [n]_K \\ m < n }} \mu(f^{n-m} \otimes f) 
&= \sum_{n=1}^{K}\sum_{m=0}^{n-1} \mu(f^{n-m} \otimes f) +  \sum_{n=K+1}^{N-1}\sum_{m=n-K}^{n-1} \mu(f^{n-m} \otimes f) \\
&= \sum_{n=1}^{K}\sum_{i=1}^{n} \mu(f^{i} \otimes f) +  \sum_{n=K+1}^{N-1}\sum_{i=1}^{K} \mu(f^{i} \otimes f) \\
&= \sum_{i=1}^K (K+1-i)\mu(f^{i} \otimes f) + (N-1-K) \sum_{i=1}^K \mu(f^{i} \otimes f) \\
&=  \sum_{i=1}^K (N-i)\mu(f^i \otimes f).
\end{align*}
In a similar fashion,
\begin{align*}
\sum_{n=0}^{N-2} \sum_{\substack{m \in [n]_K \\ m > n }} \mu(f \otimes f^{m-n}) 
& = \sum_{n=0}^{N-2-K} \sum_{m = n+1}^{n+K} \mu(f \otimes f^{m-n}) + \sum_{n=N-1-K}^{N-2} \sum_{m = n+1}^{N-1} \mu(f \otimes f^{m-n})
\\
& = \sum_{n=0}^{N-2-K} \sum_{i = 1}^{K} \mu(f \otimes f^{i}) + \sum_{n=N-1-K}^{N-2} \sum_{i = 1}^{N-1-n} \mu(f \otimes f^{i})
\\
& = (N-1-K) \sum_{i = 1}^{K} \mu(f \otimes f^{i}) + \sum_{n=1}^{K} \sum_{i = 1}^{n} \mu(f \otimes f^{i})
\\
& = \sum_{i=1}^K (N-i)\mu(f^i \otimes f).
\end{align*}
Collecting, we arrive at the expression
\begin{align*}
\widetilde\Sigma = \mu(f \otimes f) + \frac1N\sum_{i=1}^{K} (N-i)(\mu(f \otimes f^i) + \mu(f^i \otimes f)).
\end{align*}
Hence,
\begin{align*}
\Sigma -  \widetilde\Sigma = \sum_{i= K+1}^{\infty} (\mu(f^i \otimes f) + \mu(f \otimes f^i)) + \frac1N\sum_{i= 1}^{K} i (\mu(f^i \otimes f) + \mu(f \otimes f^i)),
\end{align*}
whereupon
\begin{align*}
&\| {\tr((\Sigma -  \widetilde\Sigma )D^2A) }\|_\infty 
 \le  \|D^2A\|_\infty \sum_{1\le\alpha,\beta\le d} |( \Sigma - \widetilde\Sigma)_{\alpha\beta}| \\
&\le 2\|D^2A\|_\infty\sum_{1\le\alpha,\beta\le d}\left| \sum_{i= K+1}^{\infty}  \mu(f_{\alpha}^i f_{\beta}) + \frac1N\sum_{i= 1}^{K} i \mu(f_\alpha^if_\beta) \right| \\
&\le 2d^2C_2\Vert D^2A\Vert_\infty \left(\sum_{i=K+1}^{\infty} \rho(i) +  \frac1N \sum_{i= 1}^{K} i\rho(i)\right).
\end{align*}
This finishes the proof of the lemma.
\end{proof}

\begin{lem}\label{lem4}
The expression in~\eqref{eq4} satisfies the bound
\begin{align*}
& \mu\!\left(\left| \tr \! \left(\left(  \widetilde\Sigma - \frac{1}{N}\sum_{n=0}^{N-1}\sum_{m \in [n]_K} f^n \otimes f^m \right)\! D^2A(W) \right)\right|\right)
\\
\le\ &11d^2\max\{C_2,\sqrt{C_4}\}\| D^2 A\|_{\infty}\frac{\sqrt{K+1}}{\sqrt{N}}\sqrt{\sum_{i=0}^{N-1} (i+1)\rho(i)}.
\end{align*}

\end{lem}

Observe that Proposition~\ref{prop3} follows immediately from Lemmas~\ref{lem5} and~\ref{lem4}. It remains to prove the latter.
\begin{proof}[Proof of Lemma~\ref{lem4}]
We begin by estimating
\begin{align*}
&  \mu\!\left(\left| \tr \! \left(\left(  \widetilde\Sigma - \frac{1}{N}\sum_{n=0}^{N-1}\sum_{m \in [n]_K} f^n \otimes f^m \right)\! D^2A(W) \right)\right|\right)
\\
& \le \| D^2A\|_{\infty} \sum_{1\le \alpha,\beta\le d} \mu \!\left(\left|\left( \widetilde\Sigma - \frac{1}{N}\sum_{n=0}^{N-1}\sum_{m \in [n]_K} f^n \otimes f^m \right)_{\alpha\beta}\right|\right)
\\
& =
\| D^2A\|_{\infty} \sum_{1\le \alpha,\beta\le d} \frac1N\, \mu\!\left( \left|\mu \!\left(\sum_{n=0}^{N-1}\sum_{m \in [n]_K} f^n_{\alpha} f^m_{\beta}\right) - \sum_{n=0}^{N-1}\sum_{m \in [n]_K} f^n_{\alpha} f^m_{\beta} \right|\right)
\\
& = 
 \| D^2A\|_{\infty} \sum_{1\le \alpha,\beta\le d} \frac1N\, \mu\!\left( \left|\sum_{n=0}^{N-1}\sum_{m \in [n]_K} \left(f^n_{\alpha} f^m_{\beta} - \mu(f^n_{\alpha} f^m_{\beta})\right) \right|\right)
 \\
& \le 
 \| D^2A\|_{\infty} \sum_{1\le \alpha,\beta\le d} \frac1N\, \sqrt{\mu\!\left( \left(\sum_{n=0}^{N-1}\sum_{m \in [n]_K} \left(f^n_{\alpha} f^m_{\beta} - \mu(f^n_{\alpha} f^m_{\beta})\right) \right)^2\right)}
 \\
& = 
 \| D^2A\|_{\infty} \sum_{1\le \alpha,\beta\le d} \frac1N\, \sqrt{\sum_{n=0}^{N-1}\sum_{m\in [n]_{K}}\sum_{k=0}^{N-1}\sum_{l\in [k]_{K}} (\mu\!\left(f_{\alpha}^{n}f_{\beta}^{m}f_{\alpha}^{k}f_{\beta}^{l} \right)\!-\mu\!\left(f_{\alpha}^{n}f_{\beta}^{m}\right)\!\mu\!\left(f_{\alpha}^{k}f_{\beta}^{l}\right))},
\end{align*}
where the second last line uses Jensen's inequality.
Thus,
\beq
\eqref{eq4}\le  \| D^2A\|_{\infty} \sum_{1\le \alpha,\beta\le d}\frac{1}{N}\,\sqrt{\sum_{n=0}^{N-1}\sum_{m\in [n]_{K}}\sum_{k=0}^{N-1}\sum_{l\in [k]_{K}}\bigtriangleup_{nmkl}}
\label{tot error}
\eeq
where we have introduced the shorthand notation
\beqn
\bigtriangleup_{nmkl} = \bigtriangleup_{nmkl}(\alpha,\beta) = \mu (f_{\alpha}^{n}f_{\beta}^{m}f_{\alpha}^{k}f_{\beta}^{l})- \mu(f_{\alpha}^{n}f_{\beta}^{m})\mu(f_{\alpha}^{k}f_{\beta}^{l}).
\eeqn

Let us now fix $\alpha$ and $\beta$.
To facilitate bounding $\bigtriangleup_{nmkl}$, we introduce a few helpful definitions:
We say that 
\beqn
\text{$nmkl$ is a 4-index \quad if \quad $0\le n,m,k,l<N$, \quad $m\in [n]_{K}$\quad and \quad$l\in [k]_{K}$}.
\eeqn
We classify 4-indices into two cases. In Case~$1^{\circ}$, denoted $nmkl\in \mathrm{C}_{1^{\circ}}$, the numbers~$n$ and~$m$ are the smallest two or the largest two in $\{n,m,k,l \}$. That is,
\beqn
nmkl\in \mathrm{C}_{1^{\circ}} \quad\Leftrightarrow\quad \text{$\max\{n,m\} \le \min\{k,l\}$\quad or\quad $\max\{k,l\} \le \min\{n,m\}$}.
\eeqn
Case~$2^{\circ}$ is the complement, that is,
\beqn
nmkl\in \mathrm{C}_{2^{\circ}} \quad\Leftrightarrow\quad \text{$\max\{n,m\} > \min\{k,l\}$ \quad and \quad $\max\{k,l\} > \min\{n,m\}$}.
\eeqn

Sublemma \ref{cases} below gives upper bounds on $|\!\bigtriangleup_{nmkl}\!|$ separately for the two cases, using Assumption~(A1'). Let $n,m,k$ and $l$ be fixed. It will be convenient to introduce the notation 
\beqn
\{a,b,c,d\} = \{n,m,k,l\} 
\quad\text{with}\quad 
a\le b \le c \le d.
\eeqn

\begin{sublem}\label{cases}
\beqn
|\!\bigtriangleup_{nmkl}\!| \le 
\begin{cases}
(C_4+C_{2}^{2})\rho(\max \{b-a,c-b,d-c\}) & \text{in Case $1^{\circ}$}
\\
(C_4+C_{2}^{2})\rho(\max\{b-a,d-c\}) & \text{in Case $2^{\circ}$}.
\end{cases}
\eeqn
\end{sublem}
\begin{proof}
We denote both $\alpha$ and $\beta$ by $*$ to make the proof simpler to follow.

\emph{Case $1^\circ$.}  Since $\mu(f_{*}^{n}f_{*}^{m})\mu(f_{*}^{k}f_{*}^{l})=\mu(f_{*}^{a}f_{*}^{b})\mu(f_{*}^{c}f_{*}^{d})$, 
\beqn
\bigtriangleup_{nmkl} = \mu (f_{*}^{n}f_{*}^{m}f_{*}^{k}f_{*}^{l})- \mu(f_{*}^{a}f_{*}^{b})\mu(f_{*}^{c}f
_{*}^{d})  = \bigtriangleup_{abcd}.
\eeqn     
By (A1), we have 
\beqn
|\!\bigtriangleup_{abcd}\!| \le C_4\rho(c-b).
\eeqn
On the other hand, recalling that $\rho(0)= 1$ and $\rho$ is non-increasing, (A1) also yields
\beqn
\begin{split}
|\!\bigtriangleup_{abcd}\!| & \le |\mu(f_{*}^{a}f_{*}^{b}f_{*}^{c}f_{*}^{d})| + |\mu(f_{*}^{a}f_{*}^{b})\mu(f_{*}^{c}f_{*}^{d})|
\\
& \le C_4\rho(\textrm{max}\{b-a,d-c\})+ C_2\rho(b-a)C_2\rho(d-c)
\\
& \le (C_4+C_{2}^{2})\rho(\textrm{max}\{b-a,d-c\}),
\end{split}
\eeqn
as desired.

\emph{Case $2^\circ$.} Now
\beqn
\begin{split}
|\!\bigtriangleup_{nmkl}\!|&\le | \mu (f_{*}^{n}f_{*}^{m}f_{*}^{k}f_{*}^{l})| + |\mu(f_{*}^{n}f_{*}^{m})\mu(f_{*}^{k}f_{*}^{l})| 
= | \mu (f_{*}^{a}f_{*}^{b}f_{*}^{c}f_{*}^{d})| + |\mu(f_{*}^{n}f_{*}^{m})\mu(f_{*}^{k}f_{*}^{l})| 
\\
&\le C_4\rho(\max\{b-a,d-c\}) + C_2\rho(|m-n|)C_2\rho(|k-l|)
\end{split}
\eeqn
by Assumption~(A1). Note that $nmkl\in \mathrm{C}_{2^\circ}$ implies
\beqn
\max\{b-a,d-c\}\le \max\{|m-n|,|l-k|\}.  
\eeqn
Since $\rho(0)= 1$ and $\rho$ is non-increasing, it follows that 
\beqn
|\!\bigtriangleup_{nmkl}\!|\le (C_4+C_{2}^{2})\rho(\max\{b-a,d-c\}),
\eeqn
as claimed.
\end{proof}

We will separate the contributions of $4$-indices coming from Cases $1^{\circ}$ and $2^{\circ}$:
\beq\label{eq:sum_split}
\sum_{n=0}^{N-1}\sum_{m\in [n]_{K}}\sum_{k=0}^{N-1}\sum_{l\in [k]_{K}}\left|\bigtriangleup_{nmkl}\right| = \sum_{nmkl\in \mathrm{C}_{1^{\circ}}}\left|\bigtriangleup_{nmkl}\right| + \sum_{nmkl\in \mathrm{C}_{2^{\circ}}}\left|\bigtriangleup_{nmkl}\right|.
\eeq
The following two sublemmas provide upper bounds on each sum on the right side. Their proofs amount simple counting exercises and applications of Sublemma~\ref{cases}.

\begin{sublem}\label{sublem:Case1}
\beqn
\sum_{nmkl\in \mathrm{C}_{1^{\circ}}}\left|\bigtriangleup_{nmkl}\right|\le 24(C_4+C_{2}^{2})N(K+1)\sum_{L=0}^{N-1} (L+1)\rho(L).
\eeqn
\end{sublem}
\begin{proof}
We begin from the representation
\beq \label{eq:L-trick}
\sum_{nmkl\in \mathrm{C}_{1^{\circ}}}\left|\bigtriangleup_{nmkl}\right| = \sum_{L=0}^{N-1}\sum_{nmkl\in \mathrm{C}_{1^{\circ},L}}\left|\bigtriangleup_{nmkl}\right|
\eeq
where $\mathrm{C}_{1^{\circ},L}$ is the subset of $\mathrm{C}_{1^{\circ}}$ with $\max \{b-a,c-b,d-c\}=L$. For any fixed value of $L$ we see that 
$a\in \{0,...,N-1\}$, $b\in \{a,...,a+K\}\cap \{a,...,a+L\}$, $c\in \{b,...,b+L\}$ and $d\in \{c,...,c+L\}\cap \{c,...,c+K\}$. 
Furthermore, at least one of the numbers $b-a,c-b,d-c$ is exactly $L$. Thus there exist at most 
\beqn
3N(K+1)(L+1)
\eeqn
different choices for $a,b,c$ and $d$ satisfying these conditions. Since in Case $1^{\circ}$ there are at most $8$ different possible orderings of $n,m,k,l$, we deduce  
\beqn
|\mathrm{C}_{1^{\circ},L}|\le 24N(K+1)(L+1)
\eeqn
on the cardinality of $\mathrm{C}_{1^{\circ},L}$.
 Thus, by \eqref{eq:L-trick} and Sublemma \ref{cases}, we arrive at the claim.
\end{proof}

\begin{sublem}\label{sublem:Case2}
\beqn
\sum_{nmkl\in \mathrm{C}_{2^{\circ}}}\left|\bigtriangleup_{nmkl}\right|\le 32(C_4+C_{2}^{2})N(K+1)\sum_{L=0}^{N-1} (L+1)\rho(L).
\eeqn
\end{sublem}
\begin{proof}
We begin from the representation
\beq \label{eq:L-trick2}
\sum_{nmkl\in \mathrm{C}_{2^{\circ}}}\left|\bigtriangleup_{nmkl}\right| = \sum_{L=0}^{N-1}\sum_{nmkl\in \mathrm{C}_{2^{\circ},L}}\left|\bigtriangleup_{nmkl}\right|
\eeq
where $\mathrm{C}_{2^{\circ},L}$ is the subset of $\mathrm{C}_{2^{\circ}}$ with $\max \{b-a,d-c\}=L$. For any fixed value of $L$ we see that
$a\in \{0,...,N-1\}$, $b\in \{a,...,a+L\}$, $c\in \{b,...,b+K\}$, $d\in \{c,...,c+L\}$, where the condition on $c$ can be deduced from the facts that $|m-n|\le K$ and that $\{n,m\}\neq \{a,b\},\{c,d\}$. Furthermore, since either $b-a$ or $d-c$ is exactly $L$, it follows that there exist at most 
\beqn
2N(K+1)(L+1)
\eeqn
choices for $a,b,c$ and $d$ satisfying these conditions. Since in Case $2^{\circ}$ there are at most $16$ different possible orderings of $n,m,k,l$, we deduce  
\beqn
|\mathrm{C}_{2^{\circ},L}|\le 32N(K+1)(L+1)
\eeqn
on the cardinality of $\mathrm{C}_{2^{\circ},L}$. Thus, by \eqref{eq:L-trick2} and Sublemma \ref{cases}, we arrive at the claim.
\end{proof}

Combining Sublemmas~\ref{sublem:Case1} and~\ref{sublem:Case2} with~\eqref{eq:sum_split}, we arrive at the estimate
\beq\label{eq:56}
\sum_{n=0}^{N-1}\sum_{m\in [n]_{K}}\sum_{k=0}^{N-1}\sum_{l\in [k]_{K}}\left|\bigtriangleup_{nmkl}\right|
\le 56(C_4+C_{2}^{2})N(K+1)\sum_{L=0}^{N-1} (L+1)\rho(L).
\eeq

Inserting~\eqref{eq:56} into~\eqref{tot error} yields the final bound
\begin{align*}
&\mu\!\left(\left| \tr \!\left( \widetilde\Sigma - \frac{1}{N}\sum_{n=0}^{N-1}\sum_{m \in [n]_K} f^n \otimes f^m  \right)\!\right|\right)
\\
\le\ &d^2\| D^2A\|_{\infty}\frac{1}{N}\sqrt{56(C_4+C_{2}^{2})N(K+1)\sum_{L=0}^{N-1} (L+1)\rho(L)}
\\
\le\ &d^2\| D^2 A\|_{\infty}\frac{\sqrt{112\max\{C_4,C_{2}^{2}\}}\sqrt{K+1}}{\sqrt{N}}\sqrt{\sum_{L=0}^{N-1} (L+1)\rho(L)}
\\
\le\ &\sqrt{112}\,d^2\max\{C_2,\sqrt{C_4}\}\| D^2 A\|_{\infty}\frac{\sqrt{K+1}}{\sqrt{N}}\sqrt{\sum_{L=0}^{N-1} (L+1)\rho(L)}.
\end{align*}

The proof of Lemma~\ref{lem4}, hence also that of Proposition~\ref{prop3}, is now complete.
\end{proof}

\subsubsection{Finishing the proof of Theorem~\ref{thm:pre}}
In view of Proposition~\ref{prop1} on~\eqref{eq1}, Proposition~\ref{prop3} on~\eqref{eq3}, and Assumption~(A2') on~\eqref{eq2}, the proof of Theorem~\ref{thm:pre} is complete. \qed

\subsection{Proof of Theorem~\ref{thm:pre-1d}}\label{sec:proof_pre-1d}

Recall that theorem~\ref{thm:pre} concerns the expression
\beqn
\tr  \Sigma D^2A(w) - w \cdot  \nabla A(w).
\eeqn
In the special case $d=1$ this reduces to
\beqn
\sigma^2 A'(w) - wA(w),
\eeqn
by denoting $\sigma^2 = \Sigma$ and ignoring the spurious additional derivatives. We can thus reuse most parts of the proof of Theorem~\ref{thm:pre}, keeping in mind that, formally in each expression, $D^k A$ should be replaced with $A^{(k-1)}$. A small difference is that Theorem~\ref{thm:pre-1d} assumes less regularity of~$A$: this time $A$ is differentiable with $A'$ absolutely continuous. Let us now proceed to the details.

Exactly as in~\eqref{eq1}--\eqref{eq2},
\begin{align}
&|\mu(\sigma^2 A'(W) - W A(W))|\nonumber
\\
&\qquad \le \left| \mu\!\left( \frac{1}{\sqrt{N}} \sum_{n=0}^{N-1} f^n (A(W) - A(W^n) - A'(W)(W-W^n)) \right) \right| \label{eq:1d_first}\\
&\qquad\quad + \left| \mu\!\left(\left(  \sigma^2 - \frac{1}{N}\sum_{n=0}^{N-1}\sum_{m \in [n]_K} f^n f^m \right)\! A'(W) \right) \right|\label{eq:1d_second}\\
&\qquad\quad +\left| \mu\!\left(\frac{1}{\sqrt{N}}\sum_{n=0}^{N-1} f^n  A(W^n)\right) \right| \label{eq:1d_third}.
\end{align}
By Assumption~(B2'), \eqref{eq:1d_third} is bounded by $\eta(N,K)$. With minor changes, we can bound~\eqref{eq:1d_first} and~\eqref{eq:1d_second} as~\eqref{eq1} and~\eqref{eq3} were bounded in Section~\ref{sec:proof_pre}.

\begin{prop}\label{prop:1d_first}
The expression in~\eqref{eq:1d_first} is bounded by
\begin{align*}
\frac12 C_2\Vert f \Vert_{\infty} \Vert A'' \Vert_{\infty}\frac{2K+1}{\sqrt{N}} \left(\rho(0) + 2\sum_{i = 1}^{2K} \rho(i)  \right).
\end{align*}
\end{prop}
\begin{proof}
Since $A$ is differentiable and $A'$ is absolutely continuous,
\beqn
\begin{split}
A(y) & = A(x) + \int_x^{y} A'(\xi)\,d\xi = A(x) + \int_x^{y} \!\left(A'(x) + \int_{x}^\xi A''(\eta)\,d\eta\right)d\xi
\\
& = A(x) - A'(x)(x-y) + \int_x^{y} \!\! \int_{x}^\xi A''(\eta)\,d\eta\, d\xi,
\end{split}
\eeqn
so
\beqn
|A(x) - A(y) - A'(x)(x-y)| \le \left|\int_x^{y} \!\! \int_{x}^\xi A''(\eta)\,d\eta\, d\xi\right| \le \frac12\|A''\|_\infty(y-x)^2.
\eeqn
Thus, \eqref{eq:1d_first} is bounded by
\beqn
 \frac12 \|f\|_\infty\|A''\|_\infty \frac{1}{\sqrt{N}} \sum_{n=0}^{N-1}\mu((W^n-W)^2) 
\eeqn
The rest of the proof, which involves bounding $\mu((W^n-W)^2)$, is identical to that of Proposition~\ref{prop1}.
\end{proof}

\begin{prop}\label{prop:1d_second}
The expression in~\eqref{eq:1d_second} is bounded by
\beqn
\begin{split}
& 2 C_2 \Vert A' \Vert_{\infty}  \left( \sum_{i = K+1}^{\infty}  \rho(i)  + \frac1N \sum_{i = 1}^{K} i \rho(i) \right)
\\
+\ &
11\max\{C_2,\sqrt{C_4}\} \| A'\|_{\infty}\frac{\sqrt{K+1}}{\sqrt{N}}\sqrt{\sum_{i=0}^{N-1} (i+1)\rho(i)}.
\end{split}
\eeqn
\end{prop}
We omit the proof, which is identical to the proof of Proposition~\ref{prop3}.

Assumption~(B2') and Propositions~\ref{prop:1d_first} and~\ref{prop:1d_second} immediately yield the final estimate.  This finishes the proof of Theorem~\ref{thm:pre-1d}. \qed


\section{Proofs of main results in discrete time}\label{sec:proofs}

In this section we derive Theorems~\ref{thm:main} and~\ref{thm:main-1d} from the material of Section~\ref{sec:Stein} and~\ref{sec:preliminary}.

\subsection{Proof of Theorem~\ref{thm:main}}

Since the test function $h$ in the theorem is assumed three times differentiable with bounded derivatives, Lemma~\ref{lem:steinmv} shows that the function~$A$
\beqn
A(w)= -\int_0^{\infty} \! \left\{\int_{\bR^d} h(e^{-s}w + \sqrt{1 - e^{-2s}}\,z)\,\phi_\Sigma(z)\,dz - \Phi_{\Sigma}(h)\right\}  ds,
\eeqn
introduced in~\eqref{eq:A}, is a $C^3$ solution to the multivariate Stein equation~\eqref{eq:steinmv}, and
\beqn
\|D^k A\|_\infty \le k^{-1} \|D^k h\|_\infty < \infty, \quad 1\le k\le 3.
\eeqn
Note that if we can check Assumption~(A2') of Theorem~\ref{thm:pre}, all conditions of the latter theorem are then verified. To that end, we use the following observation:
\begin{lem}
Under the conditions of Theorem~\ref{thm:main},
\beqn
\left|\mu\!\left(\frac{1}{\sqrt{N}}\sum_{n=0}^{N-1} f^n \cdot  \nabla A(W^n)\right)\right|  \leq \sqrt N\tilde\rho(K).
\eeqn 
\end{lem}
\begin{proof}
Dominated convergence yields
\begin{align*}
 \nabla A(w)= -\int_0^{\infty} e^{-s} \int_{\bR^d}  \nabla h(e^{-s}w + \sqrt{1 - e^{-2s}}\,z)\,\phi_\Sigma(z)\,dz\,ds,
\end{align*}
and, by Fubini's theorem,
\begin{align*}
\mu(f^n \cdot  \nabla A(W^n))
& = 
- \int_0^{\infty} e^{-s} \int_{\bR^d} \mu(f^n \cdot  \nabla h(e^{-s}W^n + \sqrt{1 - e^{-2s}}\,z))\,\phi_\Sigma(z)\,dz\,ds.
\end{align*}
Using Assumption~(A2) of Theorem~\ref{thm:main}, we have
\beqn
\begin{split}
|\mu(f^n \cdot  \nabla A(W^n))| & \le \int_0^{\infty} e^{-s} \int_{\bR^d} |\mu(f^n \cdot  \nabla h(e^{-s}W^n + \sqrt{1 - e^{-2s}}\,z))|\,\phi_\Sigma(z)\,dz\,ds
\\
& \le \tilde\rho(K) \int_0^{\infty} e^{-s} \int_{\bR^d} \phi_\Sigma(z)\,dz\,ds = \tilde\rho(K)
\end{split}
\eeqn
for all $n$, from which the claim directly follows.
\end{proof}
Therefore, also (A2') of Theorem~\ref{thm:pre} is satisfied, taking 
\beqn
\eta(N,K) = \sqrt N\tilde\rho(K).
\eeqn
Lemma~\ref{lem:steinmv} and Theorem~\ref{thm:pre}, followed by elementary estimates, now yield 
\begin{align*}
& |\mu(h(W)) - \Phi_\Sigma(h)| 
\\
=\ & |\mu(\tr \Sigma D^2A(W) - W \cdot  \nabla A(W))| \\
\le\ & d^3C_2\Vert f \Vert_{\infty} \Vert D^3 A \Vert_{\infty}\frac{2K+1}{\sqrt{N}} \left(\rho(0) + 2\sum_{i = 1}^{2K} \rho(i)\right)
 \\
&\qquad+ 2 d^2C_2\Vert D^2 A \Vert_{\infty}  \left( \sum_{i= K+1}^{\infty}  \rho(i)  + \frac1N\sum_{i= 1}^{K} i \rho(i) \right)\\
&\qquad+ 11d^2\max\{C_2,\sqrt{C_4}\}\Vert D^2 A\Vert_{\infty}\frac{\sqrt{K+1}}{\sqrt{N}}\sqrt{\sum_{i=0}^{N-1} (i+1)\rho(i)} \\
&\qquad+ \eta(N,K)
\\
\le\ & \frac{4d^3}3 C_2\Vert f \Vert_{\infty} \Vert D^3 h \Vert_{\infty}\frac{K+1}{\sqrt{N}} \sum_{i = 0}^\infty \rho(i)
 \\
&\qquad+  d^2C_2\Vert D^2 h \Vert_{\infty}  \left( \sum_{i= K+1}^{\infty}  \rho(i)  + \frac{K+1}{\sqrt N}\sum_{i= 0}^\infty (i+1) \rho(i) \right)\\
&\qquad+ 11d^2\max\{\sqrt{C_4},C_{2}\}\Vert D^2 h\Vert_{\infty}\frac{K+1}{\sqrt{N}}\sum_{i=0}^\infty (i+1)\rho(i) \\
&\qquad+ \sqrt N\tilde\rho(K)
\\
\le\ & C_* \left(\frac{K+1}{\sqrt{N}} +    \sum_{i= K+1}^{\infty}  \rho(i) \right) + \sqrt N\tilde\rho(K)
\end{align*}
with
\beqn
C_* = 12d^3\max\{\sqrt{C_4},C_{2}\}(\Vert D^2 h\Vert_{\infty} + \Vert f \Vert_{\infty} \Vert D^3 h \Vert_{\infty})\sum_{i = 0}^\infty (i+1)\rho(i).
\eeqn
Here we used $\sum_{i=0}^\infty(i+1) \ge \rho(0) \ge 1$.
This completes the proof of Theorem~\ref{thm:main}. \qed

\subsection{Proof of Theorem~\ref{thm:main-1d}}
Recalling Lemma~\ref{lem:univ_lemma},
\beqn
d_\mathscr{W}(W,Z) \le \sup_{A\in\mathscr{F}_{\sigma^2}} |\mu(\sigma^2 A'(W) - WA(W))|,
\eeqn
where $\mathscr{F}_{\sigma^2}$ consists of those functions $A\in C^1(\bR,\bR)$ that have an absolutely continuous derivative, and satisfy the bounds
\beqn
\|A\|_\infty \le 2, 
\quad
\|A'\|_\infty \le \sqrt{2/\pi}\, \sigma^{-1}
\quad\text{and}\quad
\|A''\|_\infty \le 2 \sigma^{-2}.
\eeqn

Note that Assumption~(B2) of Theorem~\ref{thm:main-1d} concerns functions $A\in C^1(\bR,\bR)$ having an absolutely continuous derivative and satisfying $\max_{0\le k\le 2}\|A^{(k)}\|_\infty\le 1$.
\begin{lem}
Under the conditions of Theorem~\ref{thm:main-1d},
\begin{align*}
\left|\mu\!\left(\frac{1}{\sqrt{N}}\sum_{n=0}^{N-1} f^nA(W^n)\right)\right| \le 2\max\{1,\sigma^{-2}\}\sqrt N\tilde\rho(K)
\end{align*}
holds for all $A\in \mathscr{F}_{\sigma^2}$.
\end{lem}
\begin{proof}
Given $A\in \mathscr{F}_{\sigma^2}$, we have $\max_{0\le k\le 2}\|A^{(k)}\|_\infty\le 2\max\{1,\sigma^{-2}\}$. Therefore, we can apply~(B2) of Theorem~\ref{thm:main-1d} to the function $(2\max\{1,\sigma^{-2}\})^{-1}A$, which yields
\beqn
|\mu( f^n A(W^n) ) | \le 2\max\{1,\sigma^{-2}\}\tilde\rho(K)
\eeqn
for $0\le n\le N$. This proves the claim.
\end{proof}
Hence, for each $A\in\mathscr{F}_{\sigma^2}$, the assumptions of Theorem~\ref{thm:pre-1d} are satisfied; in particular, Assumption~(B2') holds with
\beqn
\eta(N,K) = 2\max\{1,\sigma^{-2}\}\sqrt N\tilde\rho(K).
\eeqn
Therefore,
\begin{align*}
d_\mathscr{W}(W,Z)
& \le \frac12 C_2\Vert f \Vert_{\infty} 2 \sigma^{-2}\frac{2K+1}{\sqrt{N}} \left(\rho(0) + 2\sum_{i = 1}^{2K} \rho(i) \right)
 \\
&\qquad\quad+ 2C_2\sqrt{2/\pi}\, \sigma^{-1}  \left( \sum_{i= K+1}^{\infty}  \rho(i)  + \frac1N\sum_{i= 1}^{K} i \rho(i) \right)\\
&\qquad\quad+ 11\max\{C_2,\sqrt{C_4}\} \sqrt{2/\pi}\, \sigma^{-1}\frac{\sqrt{K+1}}{\sqrt{N}}\sqrt{\sum_{i=0}^{N-1} (i+1)\rho(i)} \\
&\qquad\quad+ \eta(N,K)
\\
&\le  4\sigma^{-2} C_2\Vert f \Vert_{\infty} \frac{K+1}{\sqrt{N}} \sum_{i = 0}^{\infty} \rho(i)
 \\
&\qquad\quad+ 2\sigma^{-1} C_2   \left( \sum_{i= K+1}^{\infty}  \rho(i)  + \frac{K+1}{\sqrt N}\sum_{i= 1}^{\infty} (i+1) \rho(i) \right)\\
&\qquad\quad+ 9\sigma^{-1}\max\{C_2,\sqrt{C_4}\} \frac{K+1}{\sqrt{N}}\sum_{i=0}^{\infty} (i+1)\rho(i) \\
&\qquad\quad+ 2\max\{1,\sigma^{-2}\}\sqrt N\tilde\rho(K)
\\
& \le \left(4\sigma^{-2} C_2\Vert f \Vert_{\infty}  + 2\sigma^{-1} C_2 +  9\sigma^{-1}\max\{C_2,\sqrt{C_4}\}\right) \frac{K+1}{\sqrt{N}} \sum_{i = 0}^{\infty}(i+1) \rho(i)
 \\
&\qquad\quad+ 2\sigma^{-1} C_2   \sum_{i= K+1}^{\infty}  \rho(i) + 2\max\{1,\sigma^{-2}\}\sqrt N\tilde\rho(K)
\\
& \le 11\max\{\sigma^{-1},\sigma^{-2}\}\max\{C_2,\sqrt{C_4}\} \left(\Vert f \Vert_{\infty} + 1 \right) \frac{K+1}{\sqrt{N}} \sum_{i = 0}^{\infty}(i+1) \rho(i)
 \\
&\qquad\quad+ 2\sigma^{-1} C_2   \sum_{i= K+1}^{\infty}  \rho(i) + 2\max\{1,\sigma^{-2}\}\sqrt N\tilde\rho(K)
\\
& \le C_{\#}\!\left(\frac{K+1}{\sqrt{N}} + \sum_{i= K+1}^{\infty}  \rho(i) \right) + C_{\#}'\sqrt N\tilde\rho(K)
\end{align*}
where
\beqn
C_{\#} = 11 \max\{\sigma^{-1},\sigma^{-2}\} \max\{C_2,\sqrt{C_4}\} (1 + \Vert f \Vert_{\infty}) \sum_{i = 0}^{\infty}(i+1) \rho(i)
\eeqn
and
\beqn
C_{\#}' = 2\max\{1,\sigma^{-2}\}
\eeqn
Here we used $\sum_{i=0}^\infty(i+1) \ge \rho(0) \ge 1$. This completes the proof of Theorem~\ref{thm:main-1d}. \qed


\section{Proof in continuous time}\label{sec:proofs_flow}
In this section we prove Theorem~\ref{thm:main_flow}, by applying Theorem~\ref{thm:main} in the case of the time-one map~$\cT = \psi^1$.

First off, note that $\sum_{i\ge 1} i\rho(i)<\infty$ if and only if $\int_0^\infty t\rho(t)<\infty$, because $\rho$ is decreasing.

We proceed to check that Assumption~(A1) of Theorem~\ref{thm:main} follows from~(C1). This is completely elementary, the only minor remark to be made being that the function~$F^k$ depends on~$f^s$ in the interval range $s\in[k,k+1]$. This effectively implies that decorrelation between~$F$ and~$F^k$ cannot be seen unless $k>1$, which results in the rate $\rho(k-1)$ instead of $\rho$. For instance,
\beqn
|\mu(F_\alpha F_\beta)| \le \int_0^1\!\!\int_0^1 |\mu(f^s_\alpha f^t_\beta)|\,ds\,dt \le C_2\rho(0)
\eeqn
and, for $k\ge 1$,
\beqn
|\mu(F_\alpha F_\beta^k)| \le \int_0^1\!\!\int_0^1 |\mu(f^s_\alpha f^{k+t}_\beta)|\,ds\,dt \le C_2\rho(k-1).
\eeqn
Introducing
\beqn\label{eq:rho1}
\rho_1(k) =
\begin{cases}
\rho(0), & k = 0,\\
\rho(k-1), & k\ge 1,
\end{cases}
\eeqn
we have
\beqn
|\mu(F_\alpha F_\beta^k)| \le C_2\rho_1(k)
\eeqn 
for all $k\ge 0$.
For $0\le l\le m\le n$,  we likewise have
\beqn
\begin{split}
|\mu(F_\alpha F_\beta^l F_\gamma^m F_\delta^n)| 
& \le \int_0^1\!\!\int_0^1\!\!\int_0^1\!\!\int_0^1 |\mu(f^s_\alpha f^{l+t}_\beta f^{m+u}_\gamma f^{n+v}_\delta)|\,ds\,dt\,du\,dv 
\\
& \le C_4 \min\{\rho_1(l),\rho_1(n-m)\},
\end{split}
\eeqn
which is easily seen by considering separately the cases $l\ge 1$, $n-m\ge 1$, and $l = n-m = 0$.
For $0\le l\le m\le n$, we also have
\beqn
\begin{split}
& |\mu(F_\alpha F_\beta^l F_\gamma^m F_\delta^n) - \mu(F_\alpha F_\beta^l)\mu(F_\gamma^m F_\delta^n)|
\\ 
\le\ & \int_0^1\!\!\int_0^1\!\!\int_0^1\!\!\int_0^1 |\mu(f^s_\alpha f^{l+t}_\beta f^{m+u}_\gamma f^{n+v}_\delta) - \mu(f^s_\alpha f^{l+t}_\beta)\mu(f^{m+u}_\gamma f^{n+v}_\delta)|\,ds\,dt\,du\,dv 
\\
\le\ & (C_4+C_2^2) \rho_1(m-l) \le 2\max\{C_4,C_2^2\}\rho_1(m-l).
\end{split}
\eeqn
Indeed, if~$l<m$, the integrand has the upper bound~$C_4\rho_1(m-l)$, whereas in the case~$l=m$ it is bounded by~$C_4+C_2^2$.

Next, we show that Assumption~(C2) implies (A2). Recall the expression of $V^t$ from~\eqref{eq:V^t}, and that $F^k = \int_k^{k+1} f^s\,ds$. Denoting $n = \lfloor t\rfloor$, we get
\beqn
\begin{split}
V^t = \frac{1}{\sqrt N} \sum_{0\le k< n-K} F^k + \frac{1}{\sqrt N}\sum_{n+K < k < N} F^k.
\end{split}
\eeqn
Thus, $V^t = W^n$, where the right side has the expression in~\eqref{eq:W^n}, with~$\cT$ in place of~$T$, and~$F$ in place of~$f$.
In particular,
\beqn
 |\mu(F^n\cdot\nabla h(w + W^n\tau))| 
 = \left| \int_n^{n+1}\mu( f^t \cdot  \nabla h(w + V^t\tau) )\,dt \right|
\le
 \tilde\rho(K)
\eeqn
by Assumption~(C2). This verifies Assumption~(A2) of Theorem~\ref{thm:main}.

As was discussed below Theorem~\ref{thm:main_flow}, the non-degeneracy Assumption~(C3) is equivalent to~(A3), so Theorem~\ref{thm:main} applies, with $\cT$ in place of $T$, and $F$ in place of $f$. We also need to replace $\rho$ by $\rho_1$, and $C_4$ by $2\max\{C_4,C_2^2\}$ in the expression of the constant $C_*$ in~\eqref{eq:main_constant}.

Hence, we have
\beqn
\begin{split}
|\mu(h(W)) - \Phi_\Sigma(h)| 
& \le C \! \left(\frac {K+1}{\sqrt{\lfloor T\rfloor}} + \sum_{i= K+1}^{\infty}  \rho_1(i)\right)  + \sqrt{\lfloor T\rfloor} \tilde\rho(K)
\\
& \le \sqrt{2}C \! \left(\frac {K+1}{\sqrt{T}} + \sum_{i= K}^{\infty}  \rho(i)\right)  + \sqrt{T} \tilde\rho(K)
\end{split}
\eeqn
where
\beqn
\begin{split}
C & = 12d^3\max\{C_2,\sqrt{2\max\{C_4,C_2^2\}}\}\left(\Vert D^2 h\Vert_{\infty} + \Vert F \Vert_{\infty} \Vert D^3 h \Vert_{\infty} \right)\sum_{i = 0}^\infty (i+1)\rho_1(i) 
\\
& \le 12d^3 \sqrt{2}\max\{C_2,\sqrt{C_4}\}\left(\Vert D^2 h\Vert_{\infty} + \Vert f \Vert_{\infty} \Vert D^3 h \Vert_{\infty} \right) 3 \sum_{i = 0}^\infty (i+1)\rho(i)
\\
& = 3\sqrt 2 C_*,
\end{split}
\eeqn
where $C_*$ has the same exact expression as in~\eqref{eq:main_constant}.
Finally, recalling~\eqref{eq:VW}, we obtain
\beqn
|\mu(h(V)) - \Phi_\Sigma(h)| \le   6C_*\!\left(\frac {K+1}{\sqrt{T}} + \sum_{i= K}^{\infty}  \rho(i)\right)  + \sqrt{T} \tilde\rho(K) + \frac{2d \|\nabla h\|_\infty\|f\|_\infty}{\sqrt T},
\eeqn
which finishes the proof of Theorem~\ref{thm:main_flow}. \qed


\section{Verifying Assumptions~(A2) and~(B2)}\label{sec:applications}
This section is of crucial importance to our paper, without which the main theorems would be of little practical use. Working in discrete time for simplicity, we are going to develop a scheme for verifying Assumptions~(A2) and~(B2), which we will then demonstrate by examples. The other conditions of the theorems are readily satisfied by large families of dynamical systems (including our examples) and observables, so we will focus only on~(A2) and~(B2).

The discussion in Section~\ref{sec:scheme} will be somewhat informal, as the goal is to outline a plausible pathway to obtaining the type of bounds appearing in Assumptions~(A2) and~(B2). In particular, nothing will be rigorously proved there. Yet, it is an abstraction of the method used in Sections~\ref{sec:2xmod1} and~\ref{sec:billiards}, where actual proofs will be given.

\subsection{Abstract scheme}\label{sec:scheme}
Before entering the subject matter, we remind the reader that we always use the maximum over the componentwise supremum norms to measure the size of vector-valued and matrix-valued functions, etc; see the end of Section~\ref{sec:intro}.

Let us begin by recalling that, given a three times differentiable $h:\bR^d\to\bR$ with bounded derivatives, Assumption~(A2) calls for a uniform bound on
$
|\mu( f^n \cdot  \nabla h(v + W^n t) ) | 
$
for all $0\le n < N$, $0\le t\le 1$ and $v\in\bR^d$. Note that $\|D^{k}\nabla h(v+\cdot\, t)\|_\infty \le \|D^{k+1}h\|_\infty$.
On the other hand, Assumption~(B2) calls for a uniform bound on
$|\mu( f^n A(W^n) ) |$ for all $0\le n < N$ and all differentiable $A:\bR\to\bR$ with~$A'$ absolutely continuous and $\max_{0\le k\le 2}\|A^{(k)}\|_\infty \le 1$. We thus make an observation:

\bigskip
\noindent{\bf The objective.} Verifying~(A2) and~(B2) are both special cases of bounding
\beqn
|\mu(f^n \cdot B(W^n))|
\eeqn
for functions $B\in C^1(\bR^d,\bR^d)$ satisfying $\|B\|_\infty<\infty$ and $\|DB\|_\infty<\infty$.

\bigskip
To proceed, we write
\beqn
W^n = W^n_- + W^n_+
\eeqn
where
\beqn
W^n_- = \frac{1}{\sqrt{N}}\sum_{k=0}^{n-K-1} f^k
\qquad\text{and}\qquad
W^n_+ = \frac{1}{\sqrt{N}}\sum_{k=n+K+1}^{N-1} f^k.
\eeqn
Note that there are time gaps of size $K+1$ between the random variables~$W^n_-$,~$f^n$ and~$W^n_+$, as the first one depends only on $\{f^0,\dots,f^{n-K-1}\}$ and the third one only on $\{f^{n+K+1},\dots,f^{N-1}\}$. Since $\mu(f^n) = \mu(f) = 0$, we would like to conclude that
\beqn
|\mu(f^n\cdot B(W^n))| = |\mu(f^n\cdot B(W^n_- + W^n_+))|
\eeqn
is small for sufficiently large values of~$K$. We will shortly outline how this can be achieved, in several steps.

To motivate the procedure to follow, let us momentarily imagine that $W^n_-$ were constant, say identically equal to~$c\in\bR^d$. Denoting
\beqn
W^n_+ = \widetilde W^n_+ \circ T^{n+K+1}
\qquad
\text{with}
\qquad
\widetilde W^n_+ = \frac{1}{\sqrt{N}}\sum_{k=0}^{N-n-K-2} f^k,
\eeqn
we would have
\beqn
\mu(f^n\cdot B(W^n)) = \mu(f^n\cdot B(c+\widetilde W^n_+ \circ T^{n+K+1})) = \mu(f\cdot \widetilde B\circ T^{K+1})
\eeqn
where we have introduced $\widetilde B = B(c+\widetilde W^n_+)$. An appropriate correlation bound might then yield
$\lim_{K\to \infty}\mu(f\cdot \widetilde B\circ T^{K+1}) = \mu(f)\cdot\mu(\widetilde B) = 0$. Of course, there is no reason for~$W^n_-$ to be constant. The core idea to overcome this is to condition the measure $\mu$ onto subsets of~$X$ on each of which~$W^n_-$ is, at least nearly, constant. Doing so calls for a disintegration of the measure, and results in additional steps involving the conditional measures.

To describe the conditioning procedure accurately, we assume that $(X,\cB,\mu)$ is a standard probability space, and recall some basic facts from measure theory. Consider a possibly uncountable family
\beqn
\xi = \{\xi_q:q\in\cQ\}
\eeqn
of disjoint sets $\xi_q\subset X$ satisfying $\mu(X\setminus\cup_{q\in\cQ}\xi_q) = 0$. It is said that $\xi$ is a measurable partition if there exists a countable family of measurable sets $G_k\in \cB$, $k\in\bN$, such that each $\xi_q$ is of the form $\cap_{k\in\bN} \widetilde G_k$, where each $\widetilde G_k\in\{G_k,X\setminus G_k\}$. Given a measurable partition, the measure~$\mu$ admits a disintegration $\mu = \int_\cQ \nu_q\,d\lambda(q)$ into its conditional measures $\nu_q(\slot) = \mu(\slot|\,\xi_q)$ on the partition elements $\xi_q$. Here $\lambda$ is a probability measure on the index set~$\cQ$, the so-called factor measure. This allows us to write
\beqn
\int g\,d\mu = \int_\cQ\int_{\xi_q} g(x)\,d\nu_q(x)\,d\lambda(q)
\eeqn
for any integrable function~$g:X\to\bR$.

We are now prepared to present the steps of the scheme. Each step is followed by a short discussion on its implementation.

\bigskip
\noindent{\bf Step 1: Conditioning on approximate level sets of $W^n_+$.}
Given $n$, prove there exists a measurable partition $\xi=\xi(n)$ and constant vectors~$c_q\in\bR^d$ such that
\beqn
E_1 = \mu(f^n\cdot B(W^n)) - \int_\cQ\int_{\xi_q} f^n\cdot B(c_q + W^n_+)\,d\nu_q\,d\lambda(q)
\eeqn
is small.

\noindent{\it Discussion.}
This condition becomes plausible if $W^n_-$ is nearly constant on a large majority of the partition elements~$\xi_q$, say by making their diameters small.
In particular, if
\beqn
\delta = \max_{1\le\alpha\le d}\int_\cQ\int_{\xi_q} |(W^n_- - c_q)_\alpha|\,d\nu_q\,d\lambda(q)
\eeqn
is small, then
\beqn
|E_1| \le d^2\|f\|_\infty\|DB\|_\infty \delta
\eeqn
is small. Indeed,
\beqn
\begin{split}
|E_1| &\le \sum_{\alpha=1}^d \int_\cQ\int_{\xi_q} |f^n_\alpha | |B_\alpha(W^n_- + W^n_+) - B_\alpha(c_q + W^n_+)|\,d\nu_q\,d\lambda(q)
\\
&\le \sum_{\alpha=1}^d \|f_\alpha\|_\infty \int_\cQ\int_{\xi_q} \left|\int_0^1 \frac{d}{dt}B_\alpha(tW^n_- + (1-t)c_q + W^n_+)\,dt\right|d\nu_q\,d\lambda(q)
\\
&\le \sum_{\alpha=1}^d\|f_\alpha\|_\infty \int_\cQ\int_{\xi_q} \left|\sum_{\beta=1}^d\int_0^1 (W^n_- - c_q)_\alpha\,\partial_\beta B_\alpha(tW^n_- + (1-t)c_q + W^n_+)\,dt\right|d\nu_q\,d\lambda(q)
\\
&\le \sum_{\alpha=1}^d \|f_\alpha\|_\infty \sum_{\beta=1}^d \|\partial_\beta B_\alpha\|_\infty  \int_\cQ\int_{\xi_q} |(W^n_- - c_q)_\alpha|\,d\nu_q\,d\lambda(q).
\end{split}
\eeqn

\bigskip
\noindent{\bf Step 2: Memory loss for conditional measures.}
Writing
\beqn
\widetilde B_q(x) = B(c_q + \widetilde W_+^n(x)),
\eeqn
note that
\beqn
\int_\cQ\int_{\xi_q} f^n\cdot B(c_q + W^n_+)\,d\nu_q\,d\lambda(q) = \int_\cQ\int_{\xi_q} (f\cdot \widetilde B_q\circ T^{K+1})\circ T^n\,d\nu_q\,d\lambda(q).
\eeqn
Prove that 
\beqn
E_2 = \int_\cQ \! \left(\int_{\xi_q} (f\cdot \widetilde B_q\circ T^{K+1})\circ T^n\,d\nu_q - \mu( f\cdot \widetilde B_q\circ T^{K+1}) \right) d\lambda(q)
\eeqn
is small.

\noindent{\it Discussion.} This condition becomes plausible if for a large majority of $q\in\cQ$ the pushforward measure $T^n_*\nu_q$ is close to $\mu$. For instance, suppose there exist (1) 
a small $\ve\ge 0$ and a subset $\cQ_\ve\subset \cQ$ such that 
\beqn
\lambda(\cQ_\ve)\le \ve;
\eeqn
as well as (2) a small $\ve'\ge 0$ and a class of functions $\mathscr{G}\supset \{f\cdot\widetilde B_q\circ T^{K+1}:q\in\cQ\setminus\cQ_\ve\}$, such that
\beqn
\sup_{g\in\mathscr{G}}\sup_{q\in\cQ\setminus\cQ_\ve}\left|\int_{\xi_q} g\circ T^n\,d\nu_q-\int g \,d\mu\right| \le \ve'.
\eeqn
Then
\beqn
|E_2| \le 2d\|f\|_\infty\|B\|_\infty\ve + \ve'
\eeqn
is small.
Indeed,
\beqn
\begin{split}
|E_2| & \le \lambda(\cQ_\ve)\sup_{q\in\cQ_\ve} 2\|f\cdot\widetilde B_q\|_\infty + \lambda(\cQ\setminus\cQ_\ve) \sup_{g\in\mathscr{G}}\sup_{q\in\cQ\setminus\cQ_\ve}\left|\int_{\xi_q} g\circ T^n\,d\nu_q-\int g \,d\mu\right|
\\
& \le 2d\|f\|_\infty\|B\|_\infty\ve + (1-\lambda(\cQ_\ve))\ve'.
\end{split}
\eeqn

\bigskip
\noindent{\bf Step 3: Decorrelation with respect to~$\mu$.}
Prove that
\beqn
E_3 = \int_{\cQ}\mu( f\cdot \widetilde B_q\circ T^{K+1})\,d\lambda(q)
\eeqn
is small.

\noindent{\it Discussion.} We are now back in the imaginary case $W^n_- \equiv c$ discussed before Step~1. For instance, suppose there exist a small $\ve''>0$ and a class of functions $\mathscr{G}'\supset\{\widetilde B_q:q\in\cQ\}$ such that
\beqn
\sup_{g\in\mathscr{G}'}\left|f\cdot g\circ T^{K+1}\right| \le \ve''.
\eeqn
Then
\beqn
|E_3| \le \ve''
\eeqn
holds in particular.

\bigskip
\noindent{\bf Step 4: Final upper bound.} Collecting the estimates from Steps~1--3, prove that there exists a small $\tilde\rho(K)$ such that
\beqn
\max_{0\le n\le N}|\mu(f^n\cdot B(W^n))| \le \max_{0\le n\le N}(|E_1| + |E_2| + |E_3|) \le \tilde\rho(K).
\eeqn

\noindent{\it Discussion.} To get a bound that depends on~$K$ of course requires choices to be made in Steps~1--3. This is best clarified by the examples that follow.

\subsection{A gentle example: \texorpdfstring{$x\mapsto 2x\pmod 1$}{2x (mod 1)}}\label{sec:2xmod1}
The purpose of this section is to convince the reader of the practicality of the above abstract scheme by executing it in the simplest possible setting --- that of the doubling map $T:[0,1]\to[0,1]:x\mapsto 2x\pmod 1$ and a scalar-valued observable $f:[0,1]\to\bR$. Below, $\mu$ denotes the Lebesgue measure on $[0,1]$, which is invariant under~$T$.

\begin{prop}\label{prop:2xmod1}
Let~$f:[0,1]\to\bR$ be Lipschitz continuous with constant~$L$, and~$\mu(f) = 0$. Let $A:\bR\to\bR$ be a bounded differentiable function with $\|A'\|_\infty<\infty$.
Then
\begin{align}\label{eq:angled}
|\mu(f^n A(W^n))|  
\leq L\!\left(\frac{\|A\|_\infty}2 + \frac{\|A'\|_\infty\|f\|_\infty}{\sqrt{N}}\right)\!2^{-K}
\end{align}
whenever $0\le n<N$, $0\le K<N$ and $N > 0$.
\end{prop}

The following corollary is immediate:
\begin{cor}
Let~$f:[0,1]\to\bR$ be Lipschitz continuous with constant~$L$, and~$\mu(f) = 0$. Assumption~(B2) of Theorem~\ref{thm:main-1d} is satisfied with
\beqn
\tilde\rho(K) = L(\tfrac12 + \|f\|_\infty)2^{-K}.
\eeqn
\end{cor}

\begin{proof}[Proof of Proposition~\ref{prop:2xmod1}]
Let us fix $0\le n<N$. Let~$\xi$ be the measurable partition of~$[0,1]$ into the subintervals \(\xi_q = ((q-1)2^{-n}, q2^{-n}) \), $q\in \cQ = \{1,\dots,2^n\}$. The Lebesgue measure~$\mu$ admits the obvious discrete disintegration into its conditional measures on the subintervals:
\beqn
\mu = \sum_{q\in\cQ} \lambda_q\,\nu_q
\eeqn
where
\beqn
\lambda_q = \mu(\xi_q) = 2^{-n}
\quad\text{and}\quad
\nu_q(\slot) = \mu(\slot|\,\xi_q) = \frac{\mu(\slot\cap\xi_q)}{\mu(\xi_q)}.
\eeqn

We first observe that $W_-^n$ is nearly constant on each~$\xi_q$: denoting 
\beqn
c_q = \nu_q(W^n_-) = {\mu(\xi_q)}^{-1}\int_{\xi_q}W^n_-\,d\mu
\eeqn
we have
\beqn
\begin{split}
\sup_{x\in \xi_q}|W^n_-(x) - c_q| 
& \le \sup_{x,y\in \xi_q}|W^n_-(x) -W^n_-(y)| \le \frac{1}{\sqrt N}\sum_{k=0}^{n-K-1} \sup_{x,y\in \xi_q} |f(T^k x) - f(T^k y)| 
\\
& \le \frac{L}{\sqrt{N}} \sum_{k=0}^{n-K-1} 2^{k-n} 
\le \frac{L}{\sqrt{N}\,2^K}.
\end{split}
\eeqn
Consequently,
\begin{align*}
\max_{q\in\cQ}\sup_{x\in \xi_q}| A(W^n(x)) - A(c_q + W^n_+(x)) | \le   \frac{L\|A'\|_\infty}{\sqrt{N}\,2^K}.
\end{align*}
We thus have
\begin{align*}
 |\mu(f^nA(W^{n}))| 
& = \left| \sum_{q\in\cQ}\lambda_q\int_{\xi_q}  f^{n}A(W^n) \, d\nu_q \right| 
\\
& \le \left| \sum_{q\in\cQ}\lambda_q \int_{\xi_q}   f^n A(c_q+ W^n_+) \, d\nu_q \right|
+
 \frac{L\|A'\|_\infty\|f\|_\infty}{\sqrt{N}\,2^K}.
\end{align*}
Recall $W^n_+ = \widetilde W^n_+ \circ T^{n+K+1}$. Since $T^n$ maps each $\xi_q$ affinely onto $(0,1)$, it is obvious that the $n$-fold pushforward of each $\nu_q$ is the Lebesgue measure on the full interval, that is,~$T^n_*\nu_q = \mu$. Therefore,
\begin{align*}
\sum_{q\in\cQ}\lambda_q \int_{\xi_q}   f^n A(c_q+ W^n_+) \, d\nu_q 
& = \sum_{q\in\cQ} \lambda_q \int_{\xi_q}  \bigl[f A(c_q + \widetilde W^n_+ \circ T^{K+1})\bigr]\circ T^n\, d\nu_q 
\\
& = \sum_{q\in\cQ} \lambda_q\, \mu(f \widetilde A_q\circ T^{K+1}),
\end{align*}
where the functions 
\beqn
\widetilde A_q = A(c_q + \widetilde W_+^n)
\quad\text{satisfy}\quad \|\widetilde A_q\|_\infty \le \|A\|_\infty.
\eeqn
We can now appeal to the following elementary correlation bound, whose proof is standard:
\begin{lem}
If $f:X\to\bR$ is Lipschitz continuous with constant $L\ge 0$, $\mu(f) = 0$, and $\widetilde A:\bR\to\bR$ is  bounded, then
\beqn
|\mu(f \widetilde A\circ T^j)| \le \frac{L\|\widetilde A\|_\infty}{2^j}
\eeqn
holds for all $j\ge 0$.
\end{lem}
Collecting, we arrive at \eqref{eq:angled}, which finishes the proof of Proposition~\ref{prop:2xmod1}.
\end{proof}

\subsection{Dispersing billiards}\label{sec:billiards}
In this section a more interesting application of the scheme developed in Section~\ref{sec:scheme} is given. The example is of more technical nature, so we begin with an overview of its properties.
 
Consider a two-dimensional billiard on the torus with convex scatterers having~$C^3$ boundaries of strictly positive curvature. On the surface of the torus there is a point particle moving with unit speed, linearly up to elastic collisions upon meeting a scatterer. Moreover, assume that the length of the free the billiard particle between consecutive collisions is bounded both above and away from zero. The dynamics of the particle can be encoded in a map $T:X\to X$ as follows: if $x\in X$ represents the initial location of the particle at the boundary of a scatterer together with the direction of its motion --- a unit vector pointing away from the scatterer --- then~$Tx$ represents its position and direction of motion immediately after the next collision with a scatterer. The map~$T$ is invertible, and preserves a smooth probability measure~$\mu$ on~$X$. The sigma-algebra~$\cB$ is that of the Borel sets of~$X$.

The space~$X$ consists of countably many connected components: one component (actually a two-dimensional cylinder) for each scatterer, each of which is further separated into a countably infinite number of so-called ``homogeneity strips''. For a pair of points~$x,y\in X$, one defines the future separation time $s_+(x,y)$ as the smallest integer~$n\ge 0$ such that~$T^n x$ and~$T^n y$ lie in different components. The past separation time~$s_-(x,y)$ is defined similarly in terms of the inverse map~$T^{-1}$.

Through almost every point~$x\in X$ runs a local stable manifold~$W^s(x)$. This is a maximal~$C^2$ curve with the property that~$T^n W^s(x)$ is completely contained in a component of~$X$, for all~$n\ge 0$. In fact, it follows that the length of~$T^n W^s(x)$ decreases exponentially as~$n\to\infty$. Given two points~$x,y\in X$, we either have~$W^s(x)=W^s(y)$ or~$W^s(x)\cap W^s(y) = \emptyset$. The (uncountable) family of all local stable manifolds forms a measurable partition of~$X$. Local unstable manifolds have identical properties in terms of the inverse map~$T^{-1}$. 

We refer to the textbook \cite{ChernovMarkarian_2006} for more details on billiards.

\begin{defn}
A function~$g:X\to\bR$ is dynamically H\"older continuous on local unstable manifolds with base $\vartheta\in(0,1)$ and constant $H\ge 0$ if
\beqn
|g(x) - g(y)| \le H\vartheta^{s_+(x,y)}
\eeqn
whenever $x$ and $y$ belong to the same local unstable manifold. Likewise, $g$ is dynamically H\"older continuous on local stable manifolds if
\beqn
|g(x) - g(y)| \le H\vartheta^{s_-(x,y)}
\eeqn
whenever $x$ and $y$ belong to the same local stable manifold.
\end{defn}
For instance, if $g$ is H\"older continuous (in the usual sense) with exponent $r\in(0,1)$ and constant $|g|_r$, then it is dynamically H\"older continuous simultaneously on local stable and unstable manifolds with $H = |g|_r$ and $\vartheta$ determined by~$r$.

\begin{prop}\label{prop:billiards}
There exist system constants $M>0$ and $\lambda_0\in(0,1)$ such that the following holds.
Let each coordinate function of $f:X\to\bR^d$ be dynamically H\"older continuous on both local stable and unstable manifolds with the same $H\ge 0$ and $\vartheta\in(0,1)$, and $\mu(f) = 0$. Let also $B\in C^1(\bR^d,\bR^d)$ satisfy $\|B\|_\infty < \infty$ and $\|DB\|_\infty < \infty$. Then
\beqn
|\mu(f^n \cdot B(W^n))| \le M\!\left(dH\|B\|_\infty+d\|f\|_\infty\|B\|_\infty + \frac{d^2\|f\|_\infty\|DB\|_\infty}{\sqrt N} \frac{H}{1-\vartheta}\right) \max(\vartheta^{1/4},\lambda_0)^K
\eeqn
whenever $0\le n<N$, $0\le K<N$ and $N > 0$.
\end{prop}

Before proving the proposition, let us discuss its implications. First of all, the following corollary is immediate:
\begin{cor}
Let each coordinate function of $f:X\to\bR^d$ be dynamically H\"older continuous on both local stable and unstable manifolds with the same $H\ge 0$ and $\vartheta\in(0,1)$, and $\mu(f) = 0$. Assumption~(A2) of Theorem~\ref{thm:main} is satisfied with
\beqn
\tilde\rho(K) = \widetilde C\lambda^K.
\eeqn
where
\beqn
\widetilde C = M\!\left(dH\|\nabla h\|_\infty+d\|f\|_\infty\|\nabla h\|_\infty + d^2\|f\|_\infty\|D^2h\|_\infty H(1-\vartheta)^{-1}\right)
\eeqn
and
\beqn
\lambda =  \max(\vartheta^{1/4},\lambda_0).
\eeqn
\end{cor}
Secondly, Assumption~(A1) is known to be satisfied for $f$ as above; see~\cite{Stenlund_2010}, and also~\cite{ChernovMarkarian_2006}. In fact, expressions for the constants $C_2$ and $C_4$ are known, as is the fact that $\rho(i) = \lambda^i$. Thus, recalling Corollary~\ref{cor:exponential}, we obtain the next result:
\begin{thm}\label{thm:billiards}
Suppose that $f:X\to\bR^d$, in addition to being dynamically H\"older continuous as above, is not a coboundary in any direction. Then
\beqn
|\mu(h(W)) - \Phi_{\Sigma}(h)| 
\le{\rm const}\cdot \frac{\log N}{\sqrt N}
\eeqn
holds for all $N>2$, and all $h:\bR^d\to\bR$ as in Theorem~\ref{thm:main}. Furthermore, the constant on the right side equals $C_* \! \left(\frac{2}{|{\log\lambda}|} + \frac{\lambda}{\sqrt 3(1-\lambda)}\right) + \widetilde C$.
\end{thm}
Some remarks are in order. In the very interesting work~\cite{Pene_2005}, based on an adaptation of Rio's method~\cite{Rio_1996}, P\`ene has obtained the rate $O(\frac{1}{\sqrt N})$ for test functions~$h$ which are only assumed Lipschitz continuous, and observables $f$ which are H\"older continuous (in the conventional sense); see also~\cite{Pene_2002,Pene_2004} for related, earlier, results by P\`ene. In particular, the logarithmic factor in Theorem~\ref{thm:billiards} is a byproduct of Stein's method. On the other hand, Theorem~\ref{thm:billiards} to our knowledge covers the broadest class of observables. Furthermore, the constant factor in the upper bound is explicit, in that it is expressed completely in terms of system constants standard in the theory of billiards and in terms of the characteristics of $f$ and $h$. This allows to make at least two additional remarks: (1) The constant is independent of the limit covariance~$\Sigma$, and (2) inspecting the expression of $C_*$, we moreover observe that the upper bound scales as $d^3$ with increasing dimension of the vector-valued observable $f:X\to\bR^d$. We would also argue that the approach based on Stein's method is relatively simple and conceptually transparent; in particular, this ultimately allows for treating the non-stationary problem of obtaining convergence rates for compositions of changing maps instead of iterates of a fixed map. 

\medskip
We are left with the very last part of the paper:
\begin{proof}[Proof of Proposition~\ref{prop:billiards}]
Let $\xi = \{\xi_q:q\in\cQ\}$ be the measurable partition consisting of local unstable manifolds. In the case of billiards, it is technically convenient to use invariance before disintegrating the measure. Namely, for any $m\ge 0$,
\beqn
\int f^n\cdot B(W^n)\,d\mu = \int f\circ T^{n-m} \cdot B(W^n\circ T^{-m})\,d\mu =  \int_\cQ\int_{\xi_q} f\circ T^{n-m} \cdot B(W^n\circ T^{-m})\,d\nu_q\,d\lambda(q),
\eeqn
where the last expression follows from disintegrating the measure~$\mu$ with respect to~$\xi$.
Whenever $x,y\in\xi_q$, we have $s_+(T^{k-m}(x),T^{k-m}(y))\ge m-k$ for all $m\ge k$. Since the the coordinate functions $f_\alpha$ are uniformly dynamically Hölder continuous,
\beqn
|f(T^{k-m}x) - f(T^{k-m}y)| = \max_{1\le\alpha\le d} |f_\alpha(T^{k-m}x) - f_\alpha(T^{k-m}y)| \le H\vartheta^{m-k}.
\eeqn
Thus, denoting
\beqn
c_q = \int_{\xi_q} W_-^n\circ T^{-m}\,d\nu_q
\eeqn
we have
\beqn
\begin{split}
\sup_{x\in\xi_q} |W_-^n(T^{-m}x) - c_q|
& =
\sup_{x\in\xi_q} \frac{1}{\sqrt N} \sum_{k = 0}^{n-K-1} \left|f(T^{k-m}x) - \int_{\xi_q} f\circ T^{k-m}\,\nu_q\right|
\\
& \le \frac{1}{\sqrt N} \sum_{k = 0}^{n-K-1} H\vartheta^{m-k} = \frac{H}{\sqrt N} \frac{\vartheta^{m-n+K}}{\vartheta^{-1}-1}
\end{split}
\eeqn
and
\beqn
\sup_{x\in\xi_q}|B_\alpha(W^n(T^{-m}x)) - B_\alpha(c_q + W_+^n(T^{-m}x))| \le d\max_{1\le\beta\le d}\|\partial_\beta B_\alpha\|_\infty \frac{H}{\sqrt N} \frac{\vartheta^{m-n+K}}{\vartheta^{-1}-1}
\eeqn
for all $m\ge n-K$, $q\in\cQ$ and $\alpha\in\{1,\dots,d\}$. Thus, recalling~$W^n_+ = \widetilde W^n_+ \circ T^{n+K+1}$,
\beqn
\begin{split}
\int f^n\cdot B(W^n)\,d\mu & = \int_\cQ\int_{\xi_q} f\circ T^{n-m} \cdot B(W^n\circ T^{-m})\,d\nu_q\,d\lambda(q)
\\
& = \int_\cQ \int_{\xi_q} f\circ T^{n-m} \cdot B(c_q + W_+^n\circ T^{-m})\,d\nu_q\,d\lambda(q) + E_1
\\
& = \int_\cQ \int_{\xi_q} f\circ T^{n-m} \cdot B(c_q + \widetilde W_+^n\circ T^{n+K+1-m})\,d\nu_q\,d\lambda(q) + E_1
\\
& = \int_\cQ \int_{\xi_q} \bigl[ f \cdot B(c_q + \widetilde W_+^n\circ T^{K+1})\bigr]\circ T^{n-m}\,d\nu_q\,d\lambda(q) + E_1
\end{split}
\eeqn
or
\beq\label{eq:billiards_step1}
\begin{split}
\int f^n\cdot B(W^n)\,d\mu = \int_\cQ \int_{\xi_q} \bigl[ f \cdot \widetilde B_q\circ T^{K+1}\bigr]\circ T^{n-m}\,d\nu_q\,d\lambda(q) + E_1,
\end{split}
\eeq
where
\beqn
\widetilde B_q(x) = B(c_q + \widetilde W_+^n(x))
\eeqn
and
\begin{align*}
|E_1| 
& \le \sup_{q\in\cQ}\sup_{x\in\xi_q} | f(T^{n-m}x) \cdot [B(W^n(T^{-m}x)) - B(c_q + W_+^n(T^{-m}x))] |
\\
& \le \sum_{\alpha = 1}^d\sup_{q\in\cQ}\sup_{x\in\xi_q} | f_\alpha(T^{n-m}x) [B_\alpha(W^n(T^{-m}x)) - B_\alpha(c_q + W_+^n(T^{-m}x))]|
\\
& \le \|f\|_\infty\sum_{\alpha = 1}^d\sup_{q\in\cQ}\sup_{x\in\xi_q} | B_\alpha(W^n(T^{-m}x)) - B_\alpha(c_q + W_+^n(T^{-m}x))|
\\
& \le \|f\|_\infty\sum_{\alpha = 1}^d d\max_{1\le\beta\le d}\|\partial_\beta B_\alpha\|_\infty\frac{H}{\sqrt N} \frac{\vartheta^{m-n+K}}{\vartheta^{-1}-1}.
\end{align*}
The last estimate yields
\beq\label{eq:E1}
|E_1| \le \frac{d^2\|f\|_\infty\|DB\|_\infty H}{\sqrt N} \frac{\vartheta^{m-n+K}}{\vartheta^{-1}-1}.
\eeq

To estimate the first term on the right side of~\eqref{eq:billiards_step1}, we make the following observation:
\begin{lem}\label{lem:help}
The coordinate function $\widetilde B_{q\alpha}$, $\alpha\in\{1,\dots,d\}$, $q\in\cQ$, as well as the functions $f\cdot \widetilde B_q\circ T^{K+1}$, $q\in\cQ$, are uniformly dynamically H\"older continuous on local stable manifolds. More precisely,
\beqn
|\widetilde B_{q\alpha}(x) - \widetilde B_{q\alpha}(y)| \le \frac{d\|DB\|_\infty}{\sqrt N} \frac{H}{1-\vartheta}\vartheta^{s_-(x,y)}
\eeqn
and
\beqn
|(f\cdot \widetilde B_q\circ T^{K+1})(x) - (f\cdot \widetilde B_q\circ T^{K+1})(y)| \le \left(\frac{d^2H\|f\|_\infty\|DB\|_\infty}{\sqrt N} \frac{\vartheta^{K+1}}{1-\vartheta} + dH\|B\|_\infty\right)\! \vartheta^{s_-(x,y)}
\eeqn
whenever $x$ and $y$ belong to the same local stable manifold.
\end{lem}

\begin{proof}
Let $x$ and $y$ belong to the same local stable manifold. The coordinate functions of~$\widetilde B_{q}$ satisfy
\beqn
\begin{split}
|\widetilde B_{q\alpha}(x) - \widetilde B_{q\alpha}(y)| & = |B_\alpha(c_q + \widetilde W_+^n(x)) - B_\alpha(c_q + \widetilde W_+^n(y))|
\\
& \le d\|DB\|_\infty |\widetilde W_+^n(x) - \widetilde W_+^n(y)| 
\le  \frac{d\|DB\|_\infty}{\sqrt N} \sum_{k = 0}^{N-n-K-2} |f(T^kx) - f(T^ky)|
\\
& \le  \frac{d\|DB\|_\infty}{\sqrt N} \sum_{k = 0}^{N-n-K-2} H\vartheta^{s_-(T^kx,T^ky)} = \frac{d\|DB\|_\infty}{\sqrt N} \sum_{k = 0}^{N-n-K-2} H\vartheta^{s_-(x,y)+k}
\\
& \le \frac{d\|DB\|_\infty}{\sqrt N} \frac{H}{1-\vartheta}\vartheta^{s_-(x,y)},
\end{split}
\eeqn
as claimed. As an immediate consequence,
\beqn
|\widetilde B_{q\alpha}(T^{K+1} x) - \widetilde B_{q\alpha}(T^{K+1} y)| \le \frac{d\|DB\|_\infty}{\sqrt N} \frac{H\vartheta^{K+1}}{1-\vartheta}\vartheta^{s_-(x,y)}.
\eeqn
Finally,
\beqn
\begin{split}
& |(f\cdot \widetilde B_q\circ T^{K+1})(x) - (f\cdot \widetilde B_q\circ T^{K+1})(y)| 
\\
\le\ & \sum_{\alpha=1}^d \left(|f_\alpha(x)| |\widetilde B_{q\alpha}\circ T^{K+1}(x) - \widetilde B_{q\alpha}\circ T^{K+1}(y)| + |f_\alpha(x)-f_\alpha(y)||\widetilde B_{q\alpha}\circ T^{K+1}(y)|\right)
\\
\le\ & \sum_{\alpha=1}^d\left(\|f\|_\infty \frac{d\|DB\|_\infty}{\sqrt N} \frac{H\vartheta^{K+1}}{1-\vartheta} + H\|B\|_\infty\right)\! \vartheta^{s_-(x,y)},
\end{split}
\eeqn
which yields the second claim.
\end{proof}

In view of the preceding lemma, we can in~\eqref{eq:billiards_step1} take advantage of the fact that the pushforward of $\nu_q$ tends to $\mu$ in the following weak sense; see~\cite{ChernovMarkarian_2006} for the proof.

\begin{lem}\label{lem:billiards_pushforward}
There exist constants $c_0>0$, $M_0>0$ and $\vartheta_0\in(0,1)$ such that the following holds.
Suppose $g$ is dynamically H\"older continuous on local stable manifolds with base~$\vartheta$ and constant $H_g$. Then, writing $\theta_0 = \max(\vartheta_0,\vartheta^{1/2})$,
\beqn
\left|\int_{\xi_q} g\circ T^j\,d\nu_q - \int g\,d\mu\right| \le M_0(H_g+\|g\|_\infty) \theta_0^{j-c_0|{\log|\xi_q|}|}
\eeqn
holds for every~$j\ge 0$ and every~$q\in\cQ$. Here $|\xi_q|$ stands for the length of the local unstable manifold $\xi_q$.
\end{lem}
Lemmas~\ref{lem:help} and~\ref{lem:billiards_pushforward} imply
\beq\label{eq:billiards_step2}
\int_{\xi_q} \bigl[ f\cdot \widetilde B_q\circ T^{K+1}\bigr]\circ T^{n-m}\,d\nu_q = \int f\cdot \widetilde B_q\circ T^{K+1} \,d\mu + E_{2q}
\eeq
where
\beq\label{eq:E2q}
|E_{2q}| \le M_0\!\left(\frac{d^2H\|f\|_\infty\|DB\|_\infty}{\sqrt N} \frac{\vartheta^{K+1}}{1-\vartheta} + dH\|B\|_\infty+d\|f\|_\infty\|B\|_\infty\right) \! \theta_0^{n-m-c_0|{\log|\xi_q|}|}.
\eeq

We proceed to estimating the first term on the right side of~\eqref{eq:billiards_step2}, using the following correlation bound; see~\cite{ChernovMarkarian_2006,Stenlund_2010,Stenlund_2012}.

\begin{lem}\label{lem:billiards_decorrelation}
There exist constants $M_1>0$ and $\vartheta_1\in(0,1)$ such that the following holds. Suppose $g$ and $h$ are dynamically H\"older continuous on local unstable and stable manifolds, respectively, both with base $\vartheta$, and respective constants $H_g$ and $H_h$. Then, writing $\theta_1 = \max(\vartheta_1,\vartheta^{1/4})$,
\beqn
\left|\int g\,h\circ T^j\,d\mu - \int g\,d\mu\int h\,d\mu\right| \le M_1(H_g\|h\|_\infty + \|g\|_\infty H_h + \|g\|_\infty\|h\|_\infty)\theta_1^j
\eeqn
for all $j\ge 0$.
\end{lem}
Thus, recalling $\int f\,d\mu = 0$, Lemmas~\ref{lem:help} and~\ref{lem:billiards_decorrelation} imply
\beq\label{eq:billiards_step3}
\int f \cdot \widetilde B_q\circ T^{K+1} \,d\mu = E_{3q},
\eeq
where
\beq\label{eq:E3q}
|E_{3q}| \le  M_1 \! \left(dH\|B\|_\infty + \|f\|_\infty \frac{d^2\|DB\|_\infty}{\sqrt N} \frac{H}{1-\vartheta} + d\|f\|_\infty\|B\|_\infty\right) \! \theta_1^{K+1}.
\eeq

The proof is almost complete, but one issue still remains to be considered. Namely, note from~\eqref{eq:E2q} that $E_{2q}$ fails to be small if $\xi_q$ is too short. To deal with this, we need the following result; the proof can be found, e.g., in~\cite{ChernovMarkarian_2006}.
\begin{lem}
There exists a constant $M_2>0$ such that
\beqn
\int_\cQ |\xi_q|^{-1} \,d\lambda(q) \le M_2.
\eeqn
Consequently, given $\ve>0$, Markov's inequality yields
\beqn
\lambda(\{q\in\cQ\,:\, |\xi_q| \le \ve\}) \le M_2\ve.
\eeqn
\end{lem}
Denoting $\cQ_\ve = \{q\in\cQ\,:\, |\xi_q| \le \ve\}$, we thus have
\beqn
\begin{split}
\left|\int_{\cQ_\ve} \int_{\xi_q} \bigl[ f\cdot \widetilde B_q\circ T^{K+1}\bigr]\circ T^{n-m}\,d\nu_q\,d\lambda(q)\right| \le d\|f\|_\infty\|B\|_\infty M_2\ve
\end{split}
\eeqn
and, by~\eqref{eq:billiards_step2} and~\eqref{eq:billiards_step3},
\beqn
\begin{split}
\left|\int_{\cQ\setminus\cQ_\ve} \int_{\xi_q} \bigl[ f \cdot \widetilde B_q\circ T^{K+1}\bigr]\circ T^{n-m}\,d\nu_q\,d\lambda(q)\right| 
&\le  \sup_{q\in\cQ\setminus\cQ_\ve} |E_{2q}| +  \sup_{q\in\cQ\setminus\cQ_\ve} |E_{3q}|.
\end{split}
\eeqn
Recalling also~\eqref{eq:billiards_step1}, we arrive at the estimate
\beqn
\left|\int f^n\cdot B(W^n)\,d\mu\right| 
 \le |E_1| +  \sup_{q\in\cQ\setminus\cQ_\ve} |E_{2q}| +  \sup_{q\in\cQ\setminus\cQ_\ve} |E_{3q}| + d\|f\|_\infty\|B\|_\infty M_2\ve.
\eeqn
Inserting the estimates obtained in~\eqref{eq:E1}, \eqref{eq:E2q} and~\eqref{eq:E3q}, we see that $\left|\int f^n\cdot B(W^n)\,d\mu\right|$ is bounded above by
\beqn
\begin{split}
& \frac{d^2\|f\|_\infty\|DB\|_\infty H}{\sqrt N} \frac{\vartheta^{m-n+K}}{\vartheta^{-1}-1}
\\
&\qquad + M_0\!\left(\frac{d^2H\|f\|_\infty\|DB\|_\infty}{\sqrt N} \frac{\vartheta^{K+1}}{1-\vartheta} + dH\|B\|_\infty+d\|f\|_\infty\|B\|_\infty\right) \! \theta_0^{n-m-c_0|{\log\ve}|}
\\
& \qquad + M_1 \! \left(dH\|B\|_\infty + \|f\|_\infty \frac{d^2\|DB\|_\infty}{\sqrt N} \frac{H}{1-\vartheta} + d\|f\|_\infty\|B\|_\infty\right) \! \theta_1^{K+1}
\\
& \qquad + d\|f\|_\infty\|B\|_\infty M_2\ve
\\
\le\ & \frac{d^2\|f\|_\infty\|DB\|_\infty H}{\sqrt N} \frac{1}{1-\vartheta}\left(\vartheta^{m-n+K+1} + M_0 \vartheta^{K+1}\theta_0^{n-m-c_0|{\log\ve}|} + M_1\theta_1^{K+1}\right)
\\
&\qquad + \left( dH\|B\|_\infty+d\|f\|_\infty\|B\|_\infty\right)\!\left(M_0 \theta_0^{n-m-c_0|{\log\ve}|} + M_1\theta_1^{K+1} + M_2\ve\right)
\\
%
%
%
%
\le\ & \frac{d^2\|f\|_\infty\|DB\|_\infty H}{\sqrt N} \frac{1}{1-\vartheta}\left(\theta^{m-n+K+1} + M_0 \theta^{K+1+\frac12(n-m-c_0|{\log\ve}|)} + M_1\theta^{\frac14(K+1)}\right)
\\
&\qquad + \left( dH\|B\|_\infty+d\|f\|_\infty\|B\|_\infty\right)\!\left(M_0 \theta^{\frac12(n-m-c_0|{\log\ve}|)} + M_1\theta^{\frac14(K+1)} + M_2\ve\right)
\end{split}
\eeqn
for every $m\ge n-K$ and $\ve>0$,
with
$
\theta = \max(\vartheta,\vartheta_0^2,\vartheta_1^4).
$
Fixing $\ve = e^{-\tfrac 1{4c_0} K}$ and $m \in [n-\frac34K-1,n-\frac34K)$
results in
\beqn
\begin{split}
&\left|\int f^n\cdot B(W^n)\,d\mu\right|
\\
&\qquad\le \frac{d^2\|f\|_\infty\|DB\|_\infty H}{\sqrt N} \frac{1}{1-\vartheta}\left(\theta^{-\frac34 K+K} + M_0 \theta^{K+1+\frac12(\frac34K-\frac14K)} + M_1\theta^{\frac14(K+1)}\right)
\\
&\qquad\qquad + \left( dH\|B\|_\infty+d\|f\|_\infty\|B\|_\infty\right)\!\left(M_0 \theta^{\frac12(\frac34K-\frac14K)} + M_1\theta^{\frac14(K+1)} + M_2 e^{-\tfrac 1{4c_0} K}\right)
\\
&\qquad\le \frac{d^2\|f\|_\infty\|DB\|_\infty H}{\sqrt N} \frac{1}{1-\vartheta}\left(\theta^{\frac14 K} + M_0 \theta^{K+1+\frac14K} + M_1\theta^{\frac14(K+1)}\right)
\\
&\qquad\qquad + \left( dH\|B\|_\infty+d\|f\|_\infty\|B\|_\infty\right)\!\left(M_0 \theta^{\frac14 K} + M_1\theta^{\frac14(K+1)} + M_2 e^{-\tfrac 1{4c_0} K}\right).
\end{split}
\eeqn
Thus,
\beqn
\left|\int f^n\cdot B(W^n)\,d\mu\right| \le M\!\left(dH\|B\|_\infty+d\|f\|_\infty\|B\|_\infty + \frac{d^2\|f\|_\infty\|DB\|_\infty}{\sqrt N} \frac{H}{1-\vartheta}\right) \theta^{\frac14 K},
\eeqn
where
\beqn
\theta = \max(\vartheta,\vartheta_0^2,\vartheta_1^4,e^{-\tfrac 1{c_0}}).
\eeqn
and
\beqn
M = 3\max(1,M_0,M_1,M_2).
\eeqn
This proves Proposition~\ref{prop:billiards}.
\end{proof}


\appendix





\newpage

\bigskip
\bigskip
\bibliography{Stein}{}
\bibliographystyle{plainurl}


\vspace*{\fill}

\end{document}